\documentclass[twoside,a4paper,english]{article}

\usepackage{babel}

\usepackage{epsfig}
\usepackage{graphicx}
\usepackage{color}

\long\def\ammark#1{}

\usepackage{amsmath}

\usepackage{myams}

\usepackage{amstext}

\usepackage{maththm}

\usepackage{eqnthm}

\usepackage{myamsproof}

\usepackage{amssymb}
\usepackage{myfnsymbol}
\usepackage{latexsym}

\usepackage{mymath}

\def\derpar#1{\partial_{#1}}

\newcommand\lint{{\protect\mathrm{-\hskip-1.1em \int}}}

\def\wto{\rightharpoonup}

\makeatletter
\def\len{\mathop{\operator@font len}\nolimits}
\def\Len{\mathop{\operator@font Len}\nolimits}
\def\dil{\mathop{\operator@font dil}\nolimits}
\def\mod{\mathop{\operator@font mod}\nolimits}
\def\supess{\mathop{\operator@font supess}\nolimits}
\makeatother

\usepackage{varioref}

\renewcommand{\thefootnote}{\,(\arabic{footnote})}

\def\v{{v_*}}

\newcommand\deriv[2][]{\frac{d #1}{d #2}}
\newcommand\pderiv[2][]{\frac{\partial #1}{\partial #2}}

\begin{document}

{\insert\footins{\footnotesize%
\noindent \emph{Date:} \today
\emph{This work is dedicated to the memory of Anthony J. Yezzi, Sr.,
father of Anthony Yezzi, Jr., who passed away shortly after its
initial completion. May he rest in peace and be remembered always
as a loving husband, father, and grandfather who is and will continue
to be dearly missed. 
}}}

\title{Metrics in the space of curves}
\author{A.~Yezzi \thanks{Georgia Institute of Technology, Atlanta, USA}
  \and
  A.~C.~G.~Mennucci \thanks{Scuola Normale Superiore, Pisa, Italy}}
\date{}
\maketitle

\begin{abstract}
  In this paper we study geometries on the manifold of curves.

  We define a manifold $M$ where objects $c\in M$ are curves, which we
  parameterize as $c:S^1\to \real^n$ ($n\ge 2$, $S^1$ is the circle).
  Given a curve $c$, we define the tangent space $T_cM$ of $M$ at $c$
  including in it all deformations $h:S^1\to\real^n$ of $c$.

  We discuss Riemannian and Finsler metrics $F(c,h)$ on this manifold
  $M$, and in particular the case of the geometric $H^0$ metric
  $F(c,h)=\int |h|^2ds$ of normal deformations $h$ of $c$; we study
  the existence of minimal geodesics of $H^0$ under constraints; we
  moreover propose a conformal version of the $H^0$ metric.
\end{abstract}

\section{Introduction}

In this paper we study geometries on the manifold $M$ of curves.
This manifold contains curves $c$, which we parameterize as 
$c:S^1\to \real^n$ ($S^1$ is the circle).
Given a curve $c$,
we define the tangent space $T_cM$ of $M$ at $c$ including in it 
deformations $h:S^1\to\real^n$, so that an infinitesimal deformation of the
curve $c$ in direction $h$ will yield the curve 
$c(\theta)+\varepsilon h(\theta)$. 
This manifold $M$ is  the \emph{Shape Space} that is studied in
this paper.

We would like to define a \emph{Riemannian metric}
on the manifold $M$ of curves: 
this means that, given two deformations $h,k\in T_cM$, we want to define
a scalar product $\langle h,k\rangle_c$, possibly dependent on $c$.
The Riemannian metric would then entail a \emph{distance} $d(c_0,c_1)$ between
the curves in $M$, defined as the infimum of the length $\Len(\gamma)$
of all smooth paths $\gamma:[0,1]\to M$ connecting $c_0$ to $c_1$.
We call \emph{minimal geodesic} a path providing the minimum of $\Len(\gamma)$
in the class of $\gamma$ with fixed endpoints.

\smallskip

A number of methods have been proposed in Shape Analysis to
define distances between shapes, averages of shapes and optimal
morphings between shapes;
some of these approaches are reviewed in section
\S\ref{sec:approaches}.  At the same time, there has been much
previous work in Shape Optimization, for example Image Segmentation
via Active Contours, 3D Stereo Reconstruction via Deformable Surfaces;
in these later methods, many authors have defined Energy Functionals
on curves (or surfaces) and utilized the Calculus of Variations to
derive curve evolutions to minimize the Energy Functionals; often
referring to these evolutions as Gradient Flows.  For example, the
well known Geometric Heat Flow, popular for its smoothing effect on
contours, is often referred as the \emph{gradient flow for length}.

The reference to these flows as \emph{gradient flows} implies a
certain Riemannian metric on the space of curves; but this fact has
been largely overlooked.  We call this metric $H^0$ henceforth. If one
wishes to have a consistent view of the geometry of the space of
curves in both Shape Optimization and Shape Analysis, then one should
use the $H^0$ metric when computing distances, averages and morphs
between shapes.

In this paper we first introduce the metric $H^0$ in \S\ref{sec: ex
  Riem}; we immediatly remark that, surprisingly, it
does not yield a well define metric structure, since the associated
distance is identically zero%
\footnote{This striking fact was first described
  in~\cite{Mumford:Gibbs}}.  In \S\ref{sec: ana} we analyse this
metric; we show that the lower-%
semi-continuous relaxation of the associated energy functional is
identically zero (see \S\ref{sec:relaxation} and \ref{prop: relaxed 0}); 
but we prove in thm.~\ref{thm:existence} that, under
additional constraints on the curvature of admissible curves, the
metric $H^0$ admits minimum geodesics; we propose in \S\ref{pulley} an
example that justifies some of the hypotheses in \ref{thm:existence}. We
can then define in~\S\ref{sec: space of bounded curvature} the
\emph{Shape Space} of curves with bounded curvature, where the metric
$H^0$ entails a positive distance. These hypotheses on curvature,
however, are not compatible with the classical definition of a
Riemannian Geometry.

More recently, a Riemannian metric was proposed in
\cite{Michor-Mumford} for the space $M$ of curves (see
\S\ref{sec:Michor-Mumford} here); this metric may fix the above
problems; but it would significantly alter the nature of gradient
flows used thus far in various Shape Optimization problems (assuming that
one wishes to make those \emph{gradient flows} consistent with this
new metric). In this metric, distances measured between curves are
defined using first and second order derivatives of the curves (and
therefore the resulting optimality conditions involve up to fourth
order derivatives); as a consequence, flows designed to converge
towards these optimality conditions are necessarily fourth order,
thereby precluding the use of \emph{Level Set Methods} which have
become popular in the field of Computer Vision and Shape Optimization.

We propose instead in \S\ref{sec: conformal} a class of
\emph{conformal metrics} that fix the above problems while minimally
altering the earlier flows: in fact the new gradient flows will amount
to a simple time reparameterization of the earlier flows.  In addition
the conformal metrics that we propose have some nice numerical and
computational properties: distances measured between curves are
defined using only first order derivatives (and therefore the
resulting optimality conditions involve only second order
derivatives); as a consequence, flows designed to converge towards
these optimality conditions are second order, thereby allowing the use
of Level Set methods: we indeed show such an implementation and a
numerical example in \S\ref{sec: numerical}.  We also proposed in
\cite{YM:eusipco} a differential operator that is well adapted to the
problem at hand: we review it here as well, in \S\ref{sec: geo par}.

\subsection{Riemannian and Finsler geometries}

Let  $M$ be a smooth connected differentiable manifold.%
\footnote{If $M$ is infinite dimensional, then we suppose that 
  it is modeled on a Banach or Hilbert separable space
  (see \cite{Lang:FDG}, ch.II).}
For any $c\in M$,
let $T_cM$ be \emph{the tangent space at $c$}, that is
the set of all vectors tangent to $M$ at $c$;
let $TM$ be the tangent bundle, that is the bundle of all tangent spaces.

\begin{Definition}\label{def: norm}
  Let $X$ be a vector space; a norm $|\cdot|$ satisfies
  \begin{enumerate}
  \item $|v|$ is positive homogenous, i.e.
    $$ \lambda |v| = | \lambda v|~~~\forall \lambda \ge 0~~~,$$
  \item  
    \[ |v+w|\le |v|+|w|\]
    (that, by (1), is equivalent to asking that $|v|$ be convex), and
  \item $|v|=0$ only when $v=0$.
  \end{enumerate}
\end{Definition}
If the last condition is not satisfied, then  $|\cdot|$ is a 
\emph{seminorm}.
\begin{Definition}\label{def: Finsler}
  We define a \emph{Finsler metric} to be a  
  function $F:TM\to\real^+$,  such that
  \begin{itemize}
  \item $F$ is lower semicontinuous and locally bounded from above,
    \ammark{che non \`e cos\`i scontato in dimensione infinita} and,
  \item 
    for all $c\in M$ , $v\mapsto F(c,v)$ is a  norm on
    $T_cM$.~%
    \footnote{Sometimes $F$ is called a ``Minkowsky norm''}
\end{itemize}

  We will also consider it to be symmetric, that is, $F(c,-v)=F(c,v)$.
\end{Definition}

We will write $|v|_c$ for $F(c,v)$, or simply $|v|$ when the base
point $c$ can be easily inferred from the context.

We define the length of any  locally Lipschitz path%
\footnote{As suggested in \cite{Michor-Mumford},
  we want to avoid referring to $\gamma$ as a \emph{curve}, because 
  confusion arises when we will introduce the manifold $M$ of 
  \emph{closed curves}.
  So  we will always talk of \emph{paths} in the infinite dimensional manifold
  $M$.
  Note also that these \emph{paths} are open-ended, 
  while \emph{curves} comprising $M$ are closed.
}
$\gamma:[0,1]\to M$ as
$$\Len \gamma =  \int_0^1 |\overdot \gamma(v)|_{\gamma(v)} ~dv$$
and the energy as
$$E(\gamma)=\int_0^1 |\overdot \gamma(v)|_{\gamma(v)}^2 ~dv$$

\smallskip 

We define the \emph{geodesic distance}
$d(x,y)$ as the infimum
\begin{equation}
d(x,y)=\inf \Len \gamma\label{eq: def d}
\end{equation}
for all locally Lipschitz paths $\gamma$ connecting $x$ to $y$.

The path connecting $x$ to $y$  that provides the $\min {E(\gamma)}$
in the class of all paths $\gamma$ connecting $x$ to $y$
is the \emph{(minimal length) geodesic connecting $x$ to $y$}.%
\footnote{As explained in the proposition \ref{prop: len and E} in
  the appendix \ref{app:finsler}, 
  we may equivalently define a minimal geodesic to be a minimum of
  $\min {\Len(\gamma)}$, that has been reparameterized to constant velocity}

Note that, in the classical books on Finsler Geometry
(see for example \cite{BaoChernShen}),
$M$ is finite dimensional, and $F(c,v)^2$ is considered to be
smooth and \emph{strongly convex} in the $v$ variable (for $v\neq 0$)
(see also \ref{sec: Geodesics and the exponential map});
this hypothesis is not needed,
though, for the  theorems in \S\ref{app:finsler}.

\subsubsection{Riemannian geometry}

\begin{Definition}
  Suppose that $M$ is %
  modeled on a Hilbert separable space $H$.
  A \emph{Riemannian geometry} is defined by 
  associating to $M$ a
  \emph{Riemannian metric} $g$;  for any $c\in M$ , $g(c)$ is a
  positive definite
  bilinear form on the tangent space $T_cM$ of $M$ at $c$: that is, if
  $h,k$ are tangent to $M$ at $c$, then $g(c)$ defines
  a scalar product $\langle h,  k\rangle_c$.
  If the form $g$ is positive semi-definite then the
  geometry is degenerate, 
  and we will speak of a \emph{pseudo-Riemannian metric}.
  If it is positive definite, then
  tangent space $T_xM$ is isomorphic to $H$, by
  means of the metric $g$.
  See \cite{Lang:FDG}, ch.VII.
\end{Definition}
A Riemannian geometry is a special case of a Finsler geometry:
we define the norm \[|h|_c=\sqrt{\langle h, h\rangle}~~~,\]
(pseudo-Riemannian geometries produce a seminorm $|h|_c$).

In the following we often drop the base point $c$ from
$\langle h, k\rangle_c$ and $|v|_c$.

\smallskip

If $M$ is finite dimensional,  then we can write 
$$\langle h, k\rangle_c= h_i g^{i,j}(c) k_j$$
in a choice of local coordinates;
then the matrix $g^{i,j}(c)$ is smooth and positive definite.

\subsubsection{Geodesics and the exponential map}
\label{sec: Geodesics and the exponential map}
Suppose that the metric $F$ is a \emph{regular Finsler metric}, that is,
$F$ is of class $C^2$ and $F(c,\cdot)^2$ is
strongly convex (for $v\neq 0$);
such is the case when $F(x,v)^2=|v|^2=\langle v,v\rangle$ in
a smooth Riemannian manifold.

Let  
\[ \ddot \gamma = \Gamma(\dot  \gamma,\gamma) \]
be the Euler--Lagrange O.D.E.  characterizing critical paths $\gamma$ of
\[ E(\gamma)=\int_0^1 |\overdot \gamma(v)|^2 ~dv \]

Define the \emph{exponential map} $\exp_c : T_pM \to M$
as $\exp_c (\eta) = \gamma(1)$ when $\gamma$ is the
\emph{geodesic} curve solving
\begin{equation}
  \Big \{\ddot \gamma(v) = \Gamma(\dot  \gamma(v),\gamma(v)), ~~~
  \gamma(0)=c, ~~~\dot  \gamma(0)=\eta \label{eq:geodetic}
\end{equation}  
  
Then we may state this extended version of the Hopf--Rinow theorem
 
\begin{Theorem}[Hopf-Rinow]\label{Hopf-Rinow}
   Suppose that 
   $M$ is \underline{finite} dimensional and connected,
   then these are equivalent:
   \begin{itemize}
   \item 
     $(M,d)$ is complete
   \item closed bounded sets are compact
   \item 
     $M$ is \emph{geodesically complete} at a point $c$, that is,
     (\ref{eq:geodetic}) can be solved for all $v\in\real$ and all $\eta$,
     that is,       the map $\eta\mapsto \exp_c(\eta)$  is well defined 
   \item 
     for any $c$, the map $\eta\mapsto \exp_c(\eta)$  is well defined and
     surjective; 
   \end{itemize}
   and all
   those imply that $\forall x,y\in M$ there   exist a minimal geodesic 
   connecting $x$ to $y$.
 \end{Theorem}

\subsubsection{Submanifolds}

The simplest examples of Riemannian
manifolds are the submanifolds of a Hilbert space $H$.
We think of the finite dimensional case, when $H=\real^n$,
or the infinite dimensional case, when we assume that $H$ is
separable.

\begin{Proposition}\label{sotto H}
  Define the distance  $d(x,y)=\|x-y\|$ in $H$.
  Suppose that $M\subset H$ is a closed submanifold.
  We may view $M$ as a metric space,  $(M,d)$: then it is complete.
  
  We may moreover induce a Riemannian structure on $M$ using the
  scalar product of $H$: this in turn induces the geodesic distance
  $d^g$, as defined in \eqref{eq: def d}.
  Then $d^g\ge d$, and $(M,d^g)$ is complete as well.
  If $M$ is of class $C^2$,
  moreover, then $d$ and $d^g$ are locally equivalent.%
  \footnote{Proof by standard arguments, see for example 
  sec. VIII.\S6 in \cite{Lang:FDG}}
\end{Proposition}

It is not guaranteed that $d$ and $d^g$ are \emph{globally} equivalent,
as shown by this example
\begin{Example}
  \footnote{We thank A.Abbondandolo for this remark.}
  Let $H=\real^2$ and   $ M= \{ (s,\sin(s^2)) \}$.
  Let $x_n=(\sqrt{\pi n}, 0)\in M$. Then $d(x_n,x_{n+1})\to 0$ whereas
  $d^g(x_n,x_{n+1})\ge 2$.
\end{Example}

In a certain sense,  infinite dimensional Riemannian
manifolds are simpler than their corresponding
finite-dimensional: indeed, by \cite{EellsElworthy},
\begin{Theorem}[Eells-Elworthy]\label{EellsElworthy}
  Any smooth differentiable manifold $M$ modeled on an 
  infinite dimensional Hilbert space $H$
  may be embedded as an open subset of a Hilbert space.
\end{Theorem}

With respect to geodesics, the matter is though much more complicated.
Suppose $M$ is \underline{infinite}
dimensional. In this case, 
if  $(M,d)$ is complete 
the  equation (\ref{eq:geodetic}) of \emph{geodesic} curves 
can be solved  for all $v\in\real$;
but (unfortunately) many other important implications
contained in the Hopf--Rinow theorem are false.

The most important example is due to Atkin
\cite{atkin75}:
\begin{Example}[Atkin] There exists an infinite dimensional 
  \underline{complete} connected
  smooth Riemannian manifold $M$ and $x,y\in M$ such that there is no geodesic
  connecting $x$ to $y$.
\end{Example}

We show a simpler example, of
an infinite dimensional Riemannian  manifold  $M$ such
that  the metric space $(M,d)$ is complete, but
there exist two points $x,y\in M$ that cannot be connected
by a minimal geodesic.
\begin{Example}[Grossman]\label{Grossman}
  \footnote{This example  is also in a remark in
    sec. VIII.\S6 in \cite{Lang:FDG}}
  Let $l^2({\natural})$ be the Hilbert space of real 
  sequences $x=(x_0,x_1,\ldots)$
  with the scalar product
  \[\langle x,y\rangle = \sum_{i=0}^\infinity x_i y_i\] 
  Let $e_i=(0,0,\ldots,0,1,0,\ldots)$ where $1$ is in the i-th position.

  \raisebox{-4em}{\includegraphics[width=0.40\textwidth]{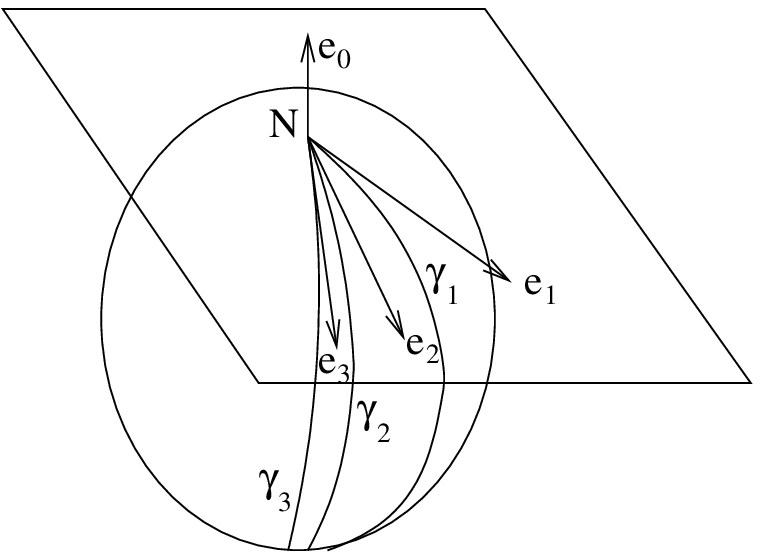}}
  \begin{minipage}[c]{0.57\textwidth}
    We build an ellipsoid 
    \[ M  = \left\{ x\in l^2 ~|~  x_0^2 + \sum _{i=1} x_i ^2 /(1+1/i)^2 = 1 \right\}\]
    in $l^2({\natural})$.
    
    Since $M$ is closed, then it is complete (with the induced metric).
  \end{minipage}

  Let    $N=e_0=(1,0,0,\ldots)$,  $S=-e_0=(-1,0,0,\ldots)$.

  Let $\gamma_i$ be the geodesic starting from
  $\gamma_i(0)=N$ and with starting speed $\dot
  \gamma_i(0)= e_i$; then there exists a first $\lambda_i>0$
  such that $\gamma_i(\lambda_i)=S$
  (moreover $\mbox{Len}(\gamma_i)=\lambda_i$).
  
  Then
  $\mbox{Len}(\gamma_i) \to \pi$,  but the sequence $\gamma_i$
  does not have a limit.

  Note that we may think of using weak convergence: but 
  $\gamma_i$ weakly converges  to the diameter;  and $e_i$
  weakly converges  to $0$.
\end{Example}

See also~\cite{ekeland78:hopf_rinow}.

\smallskip

It is then, in general, quite difficult to prove that 
an infinite dimensional manifold admits minimal geodesics
(even when it is known to be complete); a known result is

\begin{Theorem}[Cartan-Hadamard]
  Suppose that 
  $M$ is  connected, simply connected and has seminegative curvature;
  then these are equivalent:
  \begin{itemize}
  \item 
    $(M,d)$ is complete
  \item 
    for a $c$, the map $\eta\to \exp_c(\eta)$  is well defined 
  \end{itemize}
  and then there exists an unique minimal geodesic connecting any
  two points.
  \footnote{Corollary 3.9 and 3.11 in sec. IX.\S3 in \cite{Lang:FDG}}
\end{Theorem}

\subsection{Geometries of curves}

Now, suppose that $c(\theta)$ is an immersed curve $c:S_1\to \real^n$,
where $S^1$ is the circle; we want to define 
a geometry on $M$, the space of all such immersions $c$.

$M$ is a manifold; 
the tangent space $T_cM$ of $M$ at $c$ contains all the \emph{deformations}
$h\in T_cM$  of the curve $c$, that are all the vector fields along $c$. Then,
an infinitesimal deformation of the
curve $c$ in ``direction'' $h$ will yield (on first order) the curve 
$c(u)+\varepsilon h(u)$.

If $\gamma :[0,1]\to M$ is a path connecting curves,
then we may define a homotopy $C:S^1\times[0,1]\to\real^n$
associated to $\gamma$ by $C(\theta,v)=\gamma(v)(\theta)$ 
(more is in \S\ref{def: homotopy}).

\subsubsection{Finsler geometry of curves}
\label{sec: Finsler energy}
Any energy that we will study in this paper can be reconducted
to this general form

\begin{equation}
\label{eq: Finsler energy}
E(\gamma)=
\int_0^1 F\Big(\gamma(\cdot,v), \derpar v \gamma(\cdot,v)\Big)^2~dv
\end{equation}
where $F(c,h)$  is defined when $c$ is a curve,
and $h\in T_cM$ is a deformation of $c$; 
note that $F$ will be often a Minkowsky \emph{seminorm} and not a 
\emph{norm}\footnote{That is, it will fail to satisfy property 3
  in definition \ref{def: norm}}
on the space $M$ of immersions (see \ref{rem:seminorm}).

\medskip

We look mainly for metrics in the space $M$ that are independent
on the parameterization of the curves $c$: to this end, \cite{Michor-Mumford}
define 
\[  B_{i} = B_{i} (S^1,\real^2) = \mbox{Imm} (S^1,\real^2)/\mbox{Diff}(S^1) \]
and 
\[  B_{i,f} = B_{i,f} (S^1,\real^2) =\mbox{Imm}_f(S^1,\real^2)/\mbox{Diff}(S^1) \]
that are the quotients
of the spaces $\mbox{Imm}_f$ of smooth immersion, and of the
space $\mbox{Imm}_f (S^1,\real^2)$ of smooth \emph{free} immersion,
with $\mbox{Diff}(S^1)$ (the space of automorphisms of $S^1$).
$B_{i,f}$ is a manifold, the base of a principal fiber bundle, as
proved in \S2.4.3 in \cite{Michor-Mumford}, while $B_{i}$ is not.
Any metric that does not
depend on the parameterization of the curves $c$ (as defined in
eq.\eqref{eq: F = F phi}) may be projected to  $B_{i}$
by means of the results in \S2.5 in \cite{Michor-Mumford}
(the most important step appears also here as \ref{prop: black magic}).

\begin{Remark}[extending $M$]
  $\mbox{Imm}(S^1,M)$ is on open subset of the Banach space
  $C^1(S^1\to\real^n)$, where it is connected iff $n\ge 3$; whereas
  in the case $n=2$ of planar curves, it is divided in connected
  components each containing curves with the same 
  \emph{winding number}.

  To define a Riemannian Geometry on $M=\mbox{Imm}(S^1,M)$,
  it may be convenient to view it as a subset of an Hilbert
  space such as $H^1(S^1\to\real^n)$, and to complete it there.
\end{Remark}

\subsection{Abstract approach}
\label{sec: gen pro}

As a first part of this paper, we want to cast the problem
in an abstract setting.
There are some general properties that we may ask of a metric
defined as in sec.~\ref{sec: Finsler energy}.
We start with a fundamental property (that is a prerequisite
to most of the others).

\begin{enumerate}
\setcounter{enumi}{-1}
\item{\bf[well-posedness and existence of minimal geodesics]} The
  Finsler metric $F$ induces a well defined\footnote{That is, 
    the distance between different points is positive, and
    $d$ generates the same topology that the atlas of the manifold $M$ induces}
  geodesic distance $d$;
  $(M,d)$ is complete (or, it may be completed inside the space of
  mappings $c:S^1\to\real^n$); 
  for any two curves in $M$, there exists a minimal
  geodesic connecting them.
\end{enumerate}

Let $C$ be a minimal geodesic connecting $c_0$ and $c_1$.
We assume that $C$ is a homotopy of $c_0$ to $c_1$.

\begin{enumerate}
\item{\bf[rescaling]} \label{pro: rescaling}
  \emph{    if $\lambda>0$,
    if we rescale $c_0,c_1$ to $\lambda c_0,\lambda c_1$,
    then we would like that $\lambda C$ be a minimal geodesic.}
    
  A sufficient condition is that 
  $F(\lambda c,h)= \lambda^a F( c,h)$
  for some $a\ge 0$, and all $\lambda \ge 0$.
  \footnote{We could ask
  $F(\lambda c,h)= l(\lambda) F( c,h)$ 
  for some $l:\real^+\to\real^+$ monotone increasing,
  with $l(x)>0$ for $x>0$; but, 
  if $l$ is continuous, then $l(x)=x^a$ for some
  $a\ge 0$}

\item{\bf[euclidean invariance]}\label{pro: euclidean}
  \emph{ 
    If we apply an euclidean transformation $A$ to $c_0$, $c_1$,
    we would like that $AC$ be a minimal geodesic connecting $Ac_0$ to 
    $Ac_1$}

  If $F(A c,A  h)=F(c, h)$ for all $c,h$, then the above is satisfied.

\item{\bf[parameterization invariance]}\label{pro: invariant}
  there are two version of this: we define
  \begin{description}
  \item[curve-wise parameterization invariance]
    when the metric does not depend on the parameterization of the curve,
    that is
    \begin{equation}
      F(\tilde c, \tilde h)=F(c, h)
      \label{eq: F = F phi}
    \end{equation}
    when $\tilde c(t)=c(\varphi(t))$ and $\tilde h(t)=h(\varphi(t))$
    are reparameterizations of $c,h$
  \item[homotopy-wise parameterization invariance] 
    \emph{ Define $$\tilde C(\theta,v)= C( \varphi(\theta,v),v)$$
      where $\varphi:S^1\times[0,1] \to S^1$, $\varphi\in C^1$,
      $\varphi(\cdot,v)$ is a diffemorphism for all fixed $v$.  We
      would like that, in this case, $E(\tilde C)=E(C)$.}

    If $F$ can be written as
    \begin{equation}
      F=F(c,\pi_N h)\label{eq: F = F pi h}
    \end{equation}
    (that is, $F$ depends only on the normal part of the deformation)
    and if it satisfies \eqref{eq: F = F phi}, then, by proposition \ref{prop:
      C repar}, $E(\tilde C)=E(C)$.
  \end{description}
  In both cases, the geometric structure we are building depends
  only on the embedding of the curves, and not on the parameterization.
\label{lastprop}
\end{enumerate}

The above properties defined in 1,2 are valid for all examples that we will
show; property 3 is satisfied for some of them.
Property 0 is possibly the most important argument.

\begin{Definition}\label{def: geo ene}
  Any metric satisfying the above 1,2,3 is called a {\bf geometric
    metric}. 
\end{Definition}

Note that
\begin{Remark}\label{rem:seminorm}
  If $F$ satifies \eqref{eq: F = F pi h} then $F(c,\cdot)$ is
  necessarily a seminorm and not a norm\footnote{That is, it does not satisfy
    property 3 in the definition \ref{def: norm}} on the space $M$:
  we should then talk of a \emph{pseudo-Riemannian} geometry of curves.
  The projection of $F$ to the space $B_{i,f}$ may be nonetheless
  a norm.
\end{Remark}

\medskip

These other properties would be very important in applications in
Computer Vision.
\footnote{Some of these properties are much trickier: 
  we do not know sufficient conditions that imply them}

\begin{enumerate}
\setcounter{enumi}{3}
\item{\bf[finite projection]}
   \emph{There should be a finite dimensional approximation of our metric,
    for purposes  numerical computation.}

  A sufficient condition is that the energy $E(C)$ should be well defined and
  continuous with respect to a norm of a Sobolov space $W^{k,p}(I)$
  (with $k\in\natural$ and $p\in [1,\infinity)$):
  this would imply that we may approximate $C$ by smooth functions $C_h$
  and $E(C_h)\to E(C)$.

\item{\bf[embedding preserving]}
  \emph{ if $c_0$ and $c_1$ are embedded, we would like $C(\cdot,v)$
    to be embedded at all $v$}

\item{\bf[maximum principle]}
  In the following, suppose that curves are embedded in $\real^2$, and
  write $c\subset c'$ to mean that $c'$ is contained in the bounded
  region of plane enclosed by $c$.%
  \footnote{By Jordan's closed curve theorem, any embedded closed
  curve in the plane divides the plane in two regions, one bounded and
  one unbounded}

  \emph{ If $c_0\subset c'_0$, $c_1\subset c'_1$, 
    then we would like that exists   
    a minimal geodesic $C'$ connecting $c'_0$ to $c'_1$ such that
    $C'(\cdot,v)\subset C(\cdot,v)$ for $v\in[0,1]$.    }
    
    This is but another version of the Maximum Principle;
    it would imply that the minimal geodesic is unique. 
    It is an important prerequisite if we want to implement
    numerical algorithms by using \emph{Level Set Methods}.
    \ammark{idea: troncature: if $c'_v$ is convex, let 
      $\tilde c_v$ be the projection of $c_v$ on $c'_v$:
    if $F$ satisfies ... then the energy would decrease.}

\item{\bf{[convexity preserving]}}
  \emph{ if $c_0$ and $c_1$ are convex, we would like $C(\cdot,v)$
    to be convex at all $v$}

\item{\bf{[convex bounding]}}
  \emph{ 
    we would like that at all $v$, (the image of) the curve
  $C(\cdot,v)$
  be contained in the convex envelope of the curves $c_0,c_1$
  }

\item{\bf[translation]}
  \emph{ If $a\in\real^n$, 
    if $c_1=a+c_0$ is a translation of  $c_0$, we would like that 
   the uniform movement $C(\theta,v)=c_0(\theta)+va$ be
    a minimal geodesic from $c_0$ to $c_1$}

  \ammark{ If $c_1$ is a rotation ..... WHAT? may we reasonably expect 
  the minimal geodesic to rotate??}

\end{enumerate}

So we state the abstract problem:%
\footnote{
  Solving this problem in abstract would be comparable to what Shannon
  did for communication theory, where in \cite{Shannon:TheMatThe} he
  asserted there would exist a code for communication on a noisy
  channel, without actually showing an efficient algorithm to compute it.
}
\begin{Problem}
  Consider the space of curves $M$, and the family $\mathcal G$ of all
  Riemannian (or, regular Finsler) Geometries $F$  on $M$.%
  \footnote{$\mathcal G$ is non empty: see \cite{Lang:FDG}.
    By using \ref{EellsElworthy}, it would seem that there 
    exists Riemannian metrics $F$ on $M$ such that  $(M,F)$ is geodesically
    complete; we did not carry on a detailed proof.}

  Does there exist a metric $F\in \mathcal G$, satisfying  
  the above properties 0,1,2,3?

  Consider metrics $F\in \mathcal G$ that
  may be written in integral form
  \[ F(c,h)= \int_c f\big(c(s),d_s c(s),\ldots, 
  d_{s^j}^j c(s),\, h(s),\ldots, d_{s^i}^i h(s) \big) ds  \]
  what is the relationship between the
  degrees $i,j$ and the properties in this section?
\end{Problem}

\section{Examples, and different approaches and results}
\label{sec:approaches}
We now present some approaches and ideas that have been proposed
to define a metric and a distance on the space of curves;
we postpone exact definitions to section \S\ref{sec: notation}.

\subsection{Riemannian geometries of curves}
\label{sec: ex Riem}

A  Riemannian geometry is obtained by associating to 
$T_cM$ the scalar product of
an Hilbert space $H$ of squared integrable functions.
We actually have many choices for the definition of $H$.

\subsubsection{Parametric (non-geometric) form of $H^0$}

\begin{itemize}
\item  We may define 
  $$H=H^0(S^1\to\real^n)=L^2(S^1\to\real^n)$$ 
  endowed  with the scalar product
  \begin{equation}
    \langle h,k\rangle = \int_{S_1} \langle h(\theta),k(\theta)\rangle
    d\theta \label{eq:paramH0}  
  \end{equation}
  for all $h,k\in T_c M$;
  this is a common choice in analysis and geometry texts%
  \footnote{See ch 2.3 in \cite{WK:RieGeo}, or the long list of references at
  the end of II\S1 in \cite{Lang:FDG}};
  however,  the resulting metric is not invariant with respect to
  parameterization of curves (see (\ref{pro: invariant}) on page
  \pageref{pro: invariant}) and is therefore not geometric.

  \begin{Remark}
    This is the most common choice in numerical applications:
    each curve is numerically
    represented by a finite number $m$ of sample points; thereby
    discretizing the geometry of curves to the geometry of
    $\real^{nm}$.
  \end{Remark}

  Therefore \eqref{eq: Finsler energy} takes the form
  \begin{equation}
    \int_0^1 \|\overdot \gamma(v)\|_{H^0}^2 ~dv
    = \int_0^1\int_{S^1} | \derpar v \gamma(\theta,v)|^2 ~d\theta dv
    \label{eq: Riem energy}
  \end{equation}
  that is defined when $\gamma \in H^1([0,1]\to M)$.
  The energy of a homotopy is then
  $$E_s(C) =  \int_I |\derpar v C| ^2 $$

\ammark{
  We may call this a \emph{string energy}
  because we consider the curves $\theta\to C$ to be 
  deformations of the circle $S^1$ of mass $2\pi$; 
  the local mass is proportional to $d\theta$, and so is the energy needed
  to move a piece of string.
}

  \begin{Definition}
    We define the space $$H^{1,0}\big(I\to\real^n\big) =
    H^1\Big([0,1]\to H^0(S^1\to \real^n)\Big)$$
    and define the norm on
    the above space $H^{1,0}$ to be
    \begin{equation}
      \int_0^1 \int_{S^1}
      |\gamma(\theta,v) |^2+ |\derpar v \gamma(\theta,v) |^2~d\theta dv
      \label{eq: norm 1,0}
    \end{equation}
  \end{Definition}
  then
  \begin{Proposition}\label{prop: H 1 0 equiv homotopies}
    $H^{1,0}$ is the space of all finite energy homotopies,
    and the norm \eqref{eq: norm 1,0} above is equivalent
    to the energy \eqref{eq: Riem energy} on families $\gamma$ with
    fixed end points (by prop.~\ref{poincare}).
  \end{Proposition}

  $E_s(C)$ is strongly continuous in $H^{1,0}$ and is convex,
  and $E_s(C)$ is coercive (proof by using \ref{poincare}).
  So  $E_s$ has a very simple unique minimal geodesic,
  namely, the pointwise linear interpolation
  $$C^*(\theta,v)= (1-v) c_0(\theta)+v c_1(\theta)$$

\item
  As a second choice, we may define in $H$  the scalar product
  \begin{equation}
    \langle h,k\rangle = \int_0^l \langle h(s),k(s)\rangle  ds
    = \int_{S_1} \langle h(\theta),k(\theta)\rangle
    |\dot c(\theta) |d\theta \label{eq: scalar quasi N}
  \end{equation}
  where $l$ is the length of $c$, $ds$ is the arc infinitesimal, and
  $\dot c=\derpar\theta c$.\footnote{There is an abuse of
    notation in \eqref{eq: scalar quasi N}; we would like,
    intuitively, to define $h$ and $k$ on the \emph{immersed curve}
    $c$; to this end, we define $h(s)$ and $k(s)$ on the arc-parameterization
    $c(s):[0,l]\to\real^n$ (with $l=\len(c)$), and pull them back to
    $c(\theta):S^1\to\real^n$, where we write  $h(\theta),k(\theta)$
    instead of $h(s(\theta)),k(s(\theta))$ 
  }

  This scalar product does not depend on the parameterization of the curve $c$;
  but the resulting metric is still not invariant with respect to
  reparameterization of homotopies
  (see \ref{pro: invariant} in \pageref{pro: invariant}).

  By projecting this metric onto $B_{i,f}$ and lifting it back to 
  $M$ (using   \ref{prop: black magic}),
  we then devise an appropriate geometric metric, as follows.
\end{itemize}

\subsubsection{Geometric form of $H^0$}

  We then propose the scalar product
  \begin{equation}
    \langle h,k\rangle = \int_0^l \langle \pi_N h(s),\pi_N k(s)\rangle  ds
    \label{eq:H0}    
  \end{equation}
  where $\pi_N$ is the projection to the normal space $N$ to the curve
  (see \ref{def: T e pi T e pi N}).
  From now on, when we speak of the $H^0$ metric, we will be 
  implying this last definition.  

  Note that we may equivalently define the scalar product as in 
  \eqref{eq: scalar quasi N}, and only accept in $T_cM$ deformations
  that are orthogonal to the tangent of the curve. This would
  be akin to working in the quotient manifold $B_{i,f}$; or it may
  be viewed as a \emph{sub-Riemannian geometry} on $M$ itself.

  This geometric metric generates the energy
  \begin{equation}
    E^N(C) \defeq \int_{I}  |\pi_N \derpar v C|^2 |\dot C| \, d\theta\,dv  
    \label{eq:EN}
  \end{equation}
  which is invariant of reparameterizations of homotopies
  (see \ref{pro: invariant} in page \pageref{pro: invariant}).
  By proposition \ref{prop: black magic},
  the distance induced by this metric is equal to 
  the distance induced by the previous one~\eqref{eq: scalar quasi N}.

  Unfortunately, it has been  noted in \cite{Mumford:Gibbs}
  that the metric $H^0$ does not define a distance between curves, since
  \[\inf E^N(C) =0\]
  (see \ref{inf EN 0} and \ref{inf E 0}).
  We will study this metric further in section~\ref{sec: ana}.

\smallskip

There is a good reason to focus our attention
on the properties of this metric \eqref {eq:H0} for curves. Namely,
this is precisely the metric that is implicitly assumed in formulating
{\it gradient flows} of contour based energy functionals in the vast
literature on shape optimization. Consider for example the well known
{\it geometric heat flow} ($\partial_t C=\partial_{ss} C$)~%
\footnote{
  We will sometimes write $C_v$ for $\derpar v C$ (and so on),
  to simplify the derivations
} 
in which a curve evolves along the inward
normal with speed equal to its signed curvature. This flow is widely
considered to be the {\it gradient descent} of the Euclidean arclength
functional. Its smoothing properties have led to its widespread use
within the fields of computer vision and image processing. The only
sense, however, in which this is truly a gradient flow is with respect
to the $H^0$ metric as we see in the following calculation
(where $L(t)$ denotes the time varying arclength of an evolving curve $C(u,t)$
with parameter $u\in[0,1]$).

\begin{eqnarray*}
  L(t)  &&\hspace{-2em}
    =\int_C ds = \int_0^1 |C_{u}| \, d{u} \\
  L'(t) &&\hspace{-2em}
    = \int_0^1 \frac{C_{u t}\cdot C_{u}}{|C_{u}|} \, d{u}
    = \int_0^1 C_{tu}\cdot C_s \,d{u}
    = -\int_0^1 C_t\cdot C_{su} \,d{u}
    = -\int_C C_t\cdot C_{ss} \,ds \\
  &&\hspace{-2em}
    = -\Big\langle C_t,C_{ss} \Big\rangle_{H^0}
\end{eqnarray*}

If we were to change the metric then the inner-product shown
above would no longer correspond to the inner-product associated to
the metric. As a consequence, the above flow could no longer be
considered as the {\em gradient flow} with respect to the new metric.
In other words, the gradient flow would be different. We will consider
this consequence more at length in \S \ref{sec: conformal}.

\subsubsection{Michor-Mumford}
\label{sec:Michor-Mumford}
To overcome the pathologies of the $H^0$ metric,
Michor and Mumford \cite{Michor-Mumford} propose the metric
\begin{equation}
  G^A_c(h,k) = \int_0^1 (1+A|\kappa_c|^2)
\langle\pi_N h(u),\pi_N k(u)\rangle 
|\dot c(u) |du \label{eq:G MM}
\end{equation}
on planar curves,
where $\kappa_c(u)$ is the curvature of $c(u)$, and $A>0$ is fixed.
This may be generalized to  the energy%
\footnote{$J(C)$ is defined in \eqref{eq:J(C)};
  $H$ is the \emph{mean curvature}, defined in \ref{def di H}.
  Note that both $E^N(C)$ and $J(C)$ are invariant with respect to
  reparameterization, in the sense defined in \vref{pro: invariant}.
}
\begin{equation}\label{eq: E da MM}
  E^A(C)\defeq \int_0^1 \int_{S^1}
  |\pi_N\derpar v C| ^2 (1+A | H | ^2) |\overdot C|~d\theta dv
  =E^N(C) + A J(C)
\end{equation}

on space homotopies $C(u,v):[0,1]\times[0,1]\to\real^n$. 
This approach is discussed in \S\ref{sec: discuss MM}.

\subsubsection{Srivastava et al.}
\label{SEC:SRIV}

We consider here planar curves of length $2\pi$ and parameterized by arclength,
using the notation $\xi(s):S^1\to\real^2$. Such curves are 
Lipschitz continuous.

If $\xi\in C^1$,  then $|\dot\xi|=1$, so we may lift the equality
\begin{equation}
  \dot \xi(s) = \big(\!\cos(\theta(s)),~\sin(\theta(s))\big) \label{eq:lifting}
\end{equation}
to obtain a continuous function $\theta:\real\to\real$. 
This continuous lifting $\theta$ is unique up to addition of
a constant $2\pi h$, $h\in\integer$; and 
$\theta(s+ 2\pi) - \theta(s) = 2\pi i$, where $i\in\integer$ is the
\emph{winding number}, or \emph{rotation index} of $\xi$.
The addition of a generic constant to $\theta$ is equivalent to
a rotation of $\xi$.
We then understand that we may represent arc-parameterized curves $\xi(s)$,
up to translation, scaling, and rotations, by 
considering a suitable class of liftings $\theta(s)$ for $s\in[0,2\pi]$.

\smallskip

Two spaces are defined in \cite{Srivastava:AnalPlan};
we present the case of
\emph{``Shape Representation using Direction Functions''},
where the space of (pre)-shapes is
defined as the closed subset $M$ of $L^2=L^2([0,2\pi])$, 
\[ M=\left\{ \theta\in L^2([0,2\pi]) ~|~ \phi(\theta)=(2\pi^2,0,0) \right\} \]
where $\phi:L^2\to\real^3$ is defined by 
\[ \phi_1(\theta) = \int_0^{2\pi} \theta(s) ds  ~~,~~
   \phi_2(\theta) = \int_0^{2\pi}\!\cos\theta(s) ds ~~,~~
   \phi_3(\theta) = \int_0^{2\pi}\!\sin\theta(s) ds   \]

\smallskip

Define $Z$ as the set of representations  in $M$ of flat curves;
then  $Z$ is closed in $L^2$, and $M\setminus Z$ is a manifold:
\footnote{The details and the proof of \ref{prop: Sriva M}
  and \ref{prop: no rotation Srivastava}
  are in appendix \S\ref{sec: proofs Sriva}}
\begin{Proposition}\label{prop: Sriva M}
  By the implicit function theorem, $M\setminus Z$ is a smooth immersed
  submanifold of codimension 3 in $L^2$.
\end{Proposition}
Note that $M\setminus Z$ contains the (representation of) all smooth
immersed curves.

\smallskip

The manifold $M\setminus Z$ inherits a Riemannian structure, induced
by the scalar product of $L^2$; geodesics may be prolonged smoothly as
long as they do not meet $Z$.  Even if $M$ may not be a manifold at
$Z$, we may define the geodesic distance $d^g(x,y)$ in $M$ as the
infimum of the length of Lipschitz paths $\gamma:[0,1]\to L^2$ whose
image is contained in $M$;%
\footnote{It seems that $M$ is Lipschitz-arc-connected, so
  $d^g(x,y)<\infinity$; but we did not carry on a detailed proof} since
$d^g(x,y)\ge \|x-y\|_{L^2}$, and $M$ is closed in ${L^2}$, then the
metric space $(M,d^g)$ is complete.

We don't know if $(M,d^g)$ admits minimal geodesics, or if it
falls in the category of examples such as \ref{Grossman}.

\smallskip

For any $\theta\in M$, it is possible to reconstruct the
curve by integrating
\begin{equation}
  \xi(s) = \int_0^s \cos(\theta(t)),\sin(\theta(t))) dt \label{eq:unlifting}
\end{equation}
This means that $\theta$ identifies an unique curve (of length $2\pi$,
and arc-parameterized) up to
rotations and translations, and to the choice of the base point
$\xi(0)$; for this last reason, $M$ is called in
\cite{Srivastava:AnalPlan} a \emph{preshape space}. 
The \emph{shape space} $S$ is obtained by quotienting $M$
with the relation $\theta\sim \hat\theta$ iff
$\theta(s)=\hat \theta(s-a)+b$, $a,b\in\real$.
We do not discuss this quotient here.

\smallskip

We may represent any  Lipschitz closed arc-parameterized curve $\xi$ 
using a measurable $\theta\in M$:
{\def\arc{\mbox{arc}}
  let $\arc:S^1\to [0,2\pi)$ be the inverse of
  $\theta\mapsto(\cos(\theta),\sin(\theta))$; $\arc()$ is a Borel function;
  then $((\arc\circ\dot \xi)(s)+a)\in M$, for an $a\in\real$.}
We remark that the
measurable representation is never unique: for any measurable
$A,B\subset[0,2\pi]$ with $|A|=|B|$,
$(\theta(s)+2\pi({\mathbf 1}_A(s)-{\mathbf 1}_B(s)))$ will
as well represent $\xi$ in $M$. This implies that the family
$A_\xi$ of $\theta$ that represent the same curve $\xi$ 
is infinite. It may be then advisable to define a
quotient distance as follows:
\begin{equation}
  \label{eq:hatd}
  \hat d(\xi,\xi') \defeq \inf_{\theta\in A_\xi,
  \theta'\in A_{\xi'}} d(\theta,\theta')
\end{equation}
where $d(\theta,\theta')=\|\theta-\theta'\|_{L^2}$, or alternatively
$d=d^g$ is the geodesic distance on~$M$.

If $\xi\in C^1$, we have an unique%
\footnote{Indeed, the continuous lifting is unique up to addition of a
  constant to $\theta(s)$, which is equivalent to
  a rotation of $\xi$; and the constant is decided by 
  $\phi_1(\theta)=2\pi^2$
}  continuous representation $\theta\in M$; but note that,
even if $\xi,\xi'\in C^1$, the infimum \eqref{eq:hatd}
may not be given by the continuous representations $\theta,\theta'$
of  $\xi,\xi'$. Moreover there are curves $\xi$ that do not
admit a continuous representation $\theta$. As a consequence, it will not be
possible to define the rotation index of such curves $\xi$; indeed
we prove this result:
\begin{Proposition}\label{prop: no rotation Srivastava}
  For any $h\in\integer$, the set of closed smooth curves
  $\xi$ with \emph{rotation index} $h$, when represented in $M$
  using \eqref{eq:lifting}, is dense in $M\setminus Z$.
\end{Proposition}

\subsubsection{Higher order Riemannian geometry}

If we want an higher order model,
we may define a metric mimicking the definition of the Hilbert space
$H^1$, by defining
\begin{equation}
  \langle h,k\rangle_{H^1} = \langle h,k\rangle_{H^0}+
\langle \dot h,\dot k\rangle_{H^0} 
\label{eq:H1}
\end{equation}
We have again many different choices, since
\begin{itemize}
\item we may use in the RHS of \eqref{eq:H1}
  the parametric $H^0$ scalar product
  \eqref{eq:paramH0}, in which case the scalar product 
  $\langle h,k\rangle_{H^1}$ is the standard scalar product
  of $H^1(S^1\to\real^n)$; then
  homotopies are in the space
  \[H^1\big([0,1] \to H^1(S^1\to \real^n) \big)  \]
  with norm
  $$\int_I |\gamma|^2+
  \langle\derpar v \gamma \derpar v \gamma\rangle +
  \langle \derpar v \overdot \gamma ,  \derpar v \overdot \gamma\rangle   $$

\item we may use
  in the RHS of \eqref{eq:H1} the scalar product  \eqref{eq: scalar quasi N}
  or \eqref{eq:H0}. Unfortunately none of these choices
  is invariant with respect to reparameterization of homotopies.
\end{itemize}

We don't know of many application of this idea;
the only one may be considered to be \cite{Younes:Comp}.

\subsection{Finsler geometries of curves}
\label{sec: ex Fin}

To conclude, we present two examples of Finsler geometries of curves
that have been used (sometimes covertly) in the literature.

\subsubsection{$L^\infinity$ and  Hausdorff metric}

If we wish to define a norm on
$T_cM$ that is modeled on the norm of
the Banach space $L^\infinity(S^1\to\real^n)$,  we define
\[ F^\infinity(c,h)= \|\pi_N h\|_{L^\infinity}=\supess_\theta |\pi_N h(\theta)| \] 
This Finsler metric  is geometric.
The length of a homotopy is then
\[ \Len (C) =\int_0^1\supess_\theta |\pi_N \derpar v C(\theta,v)| dv \]

\paragraph{Hausdorff metric}

We recall the definition of the Hausdorff metric
\[ d_H(A,B) \defeq \max \left\{ \max_{x\in A} d(x,B) ,
\max_{y\in B} d(A,y)  \right\} \]
where $A,B\in\real^n$ are closed,  and
\[ d(x,A) \defeq  \min_{y\in A} |x-y| \]

Let $\Xi$ be the
collection of all compact subsets of $\real^n$.
We define the length of any continuous path
$\xi:[0,1]\to\Xi$
by using the total variation, as follows
\begin{equation}
  \label{eq: variation as len}
  \len ^H\gamma \defeq \sup_T \sum_{i=1}^j d_H\big(\xi(t_{i-1}),\xi(t_{i})\big)
\end{equation}
where the sup is carried out over all finite subsets
$T=\{t_0,\cdots,t_j\}$ of $[0,1]$ and $t_0\le\cdots\le t_j$.

The metric space $(\Xi,d_H)$ is complete and path-metric,%
\footnote{Path-metric: $d_H(A,B)=\inf \len ^H\gamma$ where the 
 infimum is computed in the class of Lip curves 
 $\gamma:[0,1]\to\Xi$ connecting $A$ to $B$
}
and it is possible to connect any two $A,B\in \Xi$
by a minimal geodesic, of length  $d_H(A,B)$.%
\footnote{To prove this, note that $(\Xi,d_H)$
  is locally compact and complete; and apply \ref{metric Hopf-Rinow}
}

Let $\Xi_c$ be the class of compact connected $A\subset\real^n$;
$\Xi_c$ is a closed subset of $(\Xi,d_H)$; $\Xi_C$ is  Lipschitz-path-connected
\footnote{That is, any $A,B\in \Xi_c$ can be connected 
by a Lipschitz arc $\gamma:[0,1]\to\Xi_c$};
for all above reasons, it is possible to connect any two $A,B\in \Xi_c$
by a minimal geodesic moving in $\Xi_c$;
but note that $\Xi_c$ is not geodesically convex in $\Xi$.%
\footnote{There exist two points $A,B\in\Xi_c$ and a minimal geodesics $\xi$ 
  connecting $A$ to $B$ in the metric space $(\Xi,d)$, such that
  the image of $\xi$ is not contained inside $\Xi_c$}
We don't know if $(\Xi_c,d_H)$ is path-metric.

\paragraph{Projection of $F^\infinity$ into Hausdorff metric space}

When we associate to a continuous curve $c\in M$ its image 
$Im(c)\subset\real^n$, we are actually defining
a natural projection
\[ Im : M \to \Xi_c  \]
This projection transforms a path $\gamma$ in $M$ 
into a path $\xi:[0,1]\to\Xi_c$; if  
the homotopy $C(\theta,v)=\gamma(v)(\theta)$
is continuous, then $\xi$ is continuous.
Moreover the projection of all embedded curves $c$ 
is dense in $\Xi_c$.

It is possible to prove that 
\[ \Len (\gamma) = \len^H(\xi) \] 
for a large class of paths; and then
the distance induced by the metric $F^\infinity$ coincides with
$d_H(Im(c_0),Im(c_1))$.

This is quite useful, since in
the metric space generated by the metric $F^\infinity$, it is possible
to find two curves that cannot be connected by a geodesic
(this is due to topological obstructions); whereas, the minimal geodesic will
exist in the space $(\Xi_c,d_H)$.

For this reason, \cite{Faugeras:App} proposed an approximation method
to compute $\len^H(\xi)$ by means of a family of energies
defined using a smooth integrand%
\footnote{The approximation is mainly
  based on the property $\|f\|_{L^p} \to_p \|f\|_{L^\infinity}$,
  for any measurable function $f$ defined on a bounded domain
}, and successively to find approximation of geodesics.

Unfortunately, the geometry in $(\Xi,d_H)$ is highly non regular:
for example, it is possible to find two compact sets such that
there are uncountably many minimal geodesics connecting them.%
\footnote{This is an unpublished result due to Alessandro Duci.}

\subsubsection{$L^1$ and  Plateau problem}
If we wish to define a geometric norm on
$T_cM$ that is modeled on the norm of
the Banach space $L^1(S^1\to\real^n)$,  we may define
the metric
\[ F^1(c,h)= \|\pi_N h\|_{L^1}=\int |\pi_N h(\theta)|
   |\dot c(\theta)|d\theta \] 

The length of a homotopy is then
\[ \Len (C) =\int_I  |\pi_N \derpar v C(\theta,v)| |\dot C(\theta,v)| 
d\theta dv\]
which coincides with
\[ \Len (C) =\int_I  | \derpar v C(\theta,v) \times 
\derpar\theta C(\theta,v)| 
d\theta dv\]
This last is easily recognizable as the surface area of the homotopy
(up to multiplicity); the problem of finding a minimal geodesic
connecting $c_0$ and $c_1$ 
in the $F^1$ metric may be reconducted to 
the Plateau problem of finding a surface 
which is an immersion of $I=S^1\times [0,1]$ and which has
fixed borders to the curves $c_0$ and $c_1$.
The Plateau problem is a wide and well studied subject
upon which Fomenko expounds in the monograph \cite{fomenko:plateau}.

\section{Basics}

\subsection{Relaxation of functionals}
\label{sec:relaxation}
To prove that
a Riemannian manifold admits minimal geodesics, we 
will study the energy $E(\gamma)$ by means of methods
in Calculus of Variations; we review some basic ideas.

\begingroup
\def\F{{\Omega}}
Let $X$ be a topological space, endowed with a topology $\tau$,
and  $\F\subset X$, and $f:\F\to\real$.

The function $f$ is \emph{lower semi continuous} if, for all $x\in X$,
$$ \liminf_{y\to x} f(y) \ge f(x)  $$

Note that, if the topology $\tau$ is not defined by a metric, 
then we can introduce a different condition:
$F$ is \emph{sequentially lower semi continuous} if, for all $x\in X$,
for all $x_n\to x$,
$$ \liminf_n f(x_n) \ge f(x)  $$

We define the \emph{sequential relaxation $\Gamma f$ of $f$ on $\F$}
  to be the function
\[\Gamma f:\overbar \F\to\real\]
that is  the supremum of all $f':\overbar \F\to\real$  
that are sequentially lower semicontinuous on 
$\overbar \F$,
and $f'\le f$ in $\F$. We have that, for all $x\in \overbar \F$,

$$\Gamma_\F f (x)= \min_{(x_n), x_n\to x} \{ \liminf_n f(x_n) \} $$
where the minimum is taken on all sequences $x_n$ converging to $x$,
with $x_n,x\in \F$.

$f$ is sequentially lower semicontinuous,  iff $\Gamma f|_{\F}=f$.

If $X$ is a metric space, then ``sequentially lower semicontinuous 
functions'' are ``lower semicontinuous 
functions'', and viceversa; so we may drop the adjective ``sequentially''.

\endgroup

\medskip

Consider again a function $f:X\to \real$: 
it is called \emph{coercive in $X$} if $\forall M\in\real$,
$$ \{ x\in X : f(x) \le M \}$$ 
is contained in a  compact set.

\begin{Proposition}\label{standard CV}
  If $f$ is coercive in $X$ and sequentally lower semi continuous
  then it admits a minimum on any closed set, that is, 
  for all $C\subset X$ closed there exists $x\in C$ s.t.
  $$f(x)=\min_{y\in C} f(y)$$
\end{Proposition}

This result is one of the pillars of the modern \emph{Calculus of
Variations}. We will see that unfortunately this result may not be applied
directly to the problem at hand (see \ref{exa: E N not coerc}).

\subsection{Curves and notations}
\label{sec: notation}

Consider in the following a curve $c(\theta)$ 
defined as $c:S^1\to \real^n$.

We write $\overdot c$ for $\derpar \theta c$.

\begin{Definition}\label{def: T e pi T e pi N}
  Suppose that $c$ is $C^1$;  or suppose that $c\in W^{1,1}_\text{loc}$,
  that  is, $c$  admits  a \emph{weak derivative} $\overdot c$.
  At all points where $\overdot c(\theta)\neq 0$, we define the
  tangent vector 
  $$T(\theta)= \frac { \overdot c (\theta)}{| \overdot c (\theta)|}$$
  At the points where $\overdot c=0$ we define $T=0$.

  Let $v\in\real^n$. We define the projection onto the normal
  space $N=T^\perp$
  $$  \pi_N v = v -\langle v , T\rangle T $$
  and on the tangent
  $$  \pi_T v = \langle v , T\rangle T $$
  so $\pi_N v+\pi_T v=v$ and $|\pi_N v|^2+|\pi_T v|^2=|v|^2$
  (that implies $|\pi_N v|^2=|v|^2- \langle v , T\rangle^2$).
\end{Definition}

If $c$ admits the \emph{weak derivative} $\derpar \theta c$
then $T$ is measurable, so $T\in L^\infinity$
and $\pi_T, \pi_N\in L^\infinity(S^1\to \real^{n\times n})$

\smallskip

A curve $c\in C^{1}(S^1\to \real^n)$
is \emph{immersed} when $|\overdot c|>0$ at  all points.
In this case, we can always define the 
arc parameter $s$ so that
\[ds= |\overdot c(\theta)| d\theta \]
and the derivation with respect to the arc parameter
\[ \derpar s  =\frac 1{|\overdot c|} \derpar \theta\]

We will also consider curves  $c\in W^{1,1}(S^1\to \real^n)$
such that $|\overdot c|>0$ at almost all points.

There are  two different definitions of \emph{curvature} of an immersed curve:
\emph{mean curvature $H$} and \emph{signed curvature $\kappa$}, which is
 defined when $c$ is valued in $\real^2$.

$H$ and $k$ are \emph{extrinsic curvatures}
\footnote{See also:
Eric W. Weisstein. "Extrinsic Curvature." From MathWorld--A
  Wolfram Web  Resource.
  \texttt{http://mathworld.wolfram.com/ExtrinsicCurvature.html} }:
they are properties of the \emph{embedding} of $c$ into $\real^n$.

\begin{Definition}[H]\label{def di H}
  If $c$ is $C^2$ regular and immersed,
  we can define the mean curvature $H$  of
  $c$ as 
  \[H=  \derpar s T =\frac 1{|\overdot c|} \derpar \theta T\]

  In general, we will say that a curve  $c\in W^{1,1}_\text{loc}(S^1)$ admits
  \emph{mean curvature in the measure sense} if there exists a 
  a vector valued  measure $H$ on $S^1$ 
  such that
    $$\int_I T(s) \derpar s\phi(s)~ds  = 
    - \int_I  \phi(s) \, H(ds)
    ~~~~\forall \phi\in C^\infinity (S^1)$$
    that is 
    $$\int_I T(\theta) \derpar \theta\phi(\theta)\,d\theta
    =     - \int_I  \phi(\theta) \,|\overdot c| H(d\theta)
    ~~~~\forall \phi\in C^\infinity (S^1)$$

    \ammark{ci va il $|\overdot c|$???}

  Note that the two definitions are related, since when $c\in C^2$,
  the measure is $H=\derpar s T\cdot ds$.
  See also \cite{simon}, \S7 and \S16.
\end{Definition}

We can then define the projection onto the curvature vector $H$ by 
\[ \pi_H v = \frac 1 {|H|^2} \langle v , H\rangle H \]

\begin{Definition}[N]\label{def di N}
  When the curve $c$ is in $\real^2$, and is immersed, 
  we can define%
  \footnote{There is a slight abuse of notation here, 
    since in the definition $N=T^\perp$ in \ref{def: T e pi T e pi N}
    we defined $N$ to be a ``vector space'' and not a ``vector''}
  a unit vector $N$ 
  such that $N\perp T$ and  $N$ is $\pi/2$ degree anticlockwise
  with respect to $T$. In this case   for any vector $V\in \real^2$,
  $\pi_N V = N \langle N,V \rangle$, and,
  \[ |V^2| = \pi_N V   \]
\end{Definition}

\begin{Definition}[$\kappa$]\label{def di k}
  if $c$ is in $\real^2$   then we can define a signed 
  scalar curvature
    \[ \kappa =\langle H , N\rangle  \]
    If $H$ is a measure, then $k$ is a real valued measure
    defined by
    \[ \kappa (A)=\sum_{i=1}^n\int_A N_i(\theta)\,  dH_i(d\theta)  ~~. \]

    Note that $|H|=|\kappa|$.

    When $c\in C^2(S^1\to\real^2)$ is immersed,
    \[\derpar s T = \kappa N ~~\text{ and  } ~~\derpar s N = -\kappa T ~~.\]
\end{Definition}

\begin{Remark}
  When the curve $c$ is in $\real^2$, and is immersed then
  \[ \langle H, v\rangle = \kappa \pi_N v \]
  whereas for immersed curves  $c$ in $\real^n$
  \begin{equation}
    |\langle H, v\rangle| \le |H| |\pi_N v| \label{eq:H,v<H,pi_nv}
  \end{equation}
  and we do not expect to have equality in general when $n\ge 3$,
  since $H$ is only a vector in the $(n-1)$-dimensional 
  subspace $N=T^\perp$.
\end{Remark}

\subsubsection{Homotopies}
\label{def: homotopy}

Let $I=S^1\times [0,1]$.
We define a homotopy to be a continuous function
$$C(\theta, v):S^1\times [0,1]\to \real ^n$$ 
This homotopy is a path connecting
$c_0=C(\cdot, 0)$ and $c_1=C(\cdot, 1)$
in the space $M$ of curves:
indeed any path $\gamma$ in $M$ is associated 
to a homotopy $C$ by  $C(\theta,v)=\gamma(v)(\theta)$.

We extend the above definitions to homotopies $C$,
isolating any curve $c$ in the homotopy by defining
$c(\theta)=C(\theta,v)$ for the corresponding fixed $v$;
for example, when
$C$ admits the \emph{weak derivative} $\derpar \theta C$ then

\[\pi_T, \pi_N\in L^\infinity(S^1\times[0,1]\to \real^{n\times n})\]

\begin{Remark}\label{extend H}
  We extend the measure $H(\cdot,v)$ on $S^1$, (that is curvature of
  any single curve $C(\cdot,v)$) to a Borel measure $\hat H$ on $I$,
  by
  \begin{equation}
   \hat H(A) = \int_0^1 H(A_v)dv\label{eq:extend H}   
 \end{equation}
  where $A_v$ is the section of $A$. 
$H$ can be defined equivalently  using the formula

    $$\int_I T(\theta,v) \derpar \theta\phi(\theta,v)  = 
    - \int_I  \phi(\theta,v) |\overdot C(\theta,v)| H(d\theta,dv)
    ~~~~\forall \phi\in C^\infinity_c (I)$$
  \end{Remark}

\noindent
Moreover we define the length
$$ \len(C)(v)=\int_{S^1}|\overdot C(\theta,v)|~d\theta $$
so that $\len(C):[0,1]\to\real^+$ is a function of $v$.

\subsection{Preliminary results}
\label{sec: para}

\begin{Definition}[Reparameterization]\label{prop: C repar}
  We define the reparameterization of $C$ to $\widetilde C$ as
  \begin{equation}
    \widetilde C(\theta,v)= C( \varphi(\theta,v),v)\label{eq: repar C}
  \end{equation}
  where $\varphi:S^1\times[0,1] \to S^1$ is $C^1$ regular
  with $\derpar \theta \varphi\neq 0$.
  Then (by direct computation)
  \begin{equation}
    \derpar \theta \widetilde C (\theta,v)=
    \derpar \theta   C( \tau,v)\derpar\theta\varphi(\theta,v)
  \end{equation}
  so that $T=\widetilde T\text{sign}(\derpar \theta \varphi) $;  and
  \begin{equation}
    \pi_{\widetilde N}\derpar v \widetilde C (\theta,v)=
    \pi_N\derpar v   C( \tau,v)
    \label{eq: pi d v C and repar}
  \end{equation}
  whereas
  \begin{equation}
    \pi_{\widetilde T}\derpar v \widetilde C (\theta,v)=
    \pi_T\derpar v   C( \tau,v)
    +\overdot C(\tau,v)\derpar v\varphi(\theta,v)
    \label{eq: pi T d v C and repar}
  \end{equation}
  where $\tau\defeq\varphi(\theta,v)$.
\end{Definition}

We may choose to reparameterize using the arclength parameter as in the
following proposition.
\begin{Proposition}[Arc parameter]\label{prop: by arc par}
  For any $C^1$ regular homotopy $C$ such that all curves $\theta\mapsto C$
  are immersed,
  there exists a $\varphi$ as in~\eqref{eq: repar C} above
  such that $|\derpar\theta\widetilde C|$ is constant in $\theta$ for any $v$
  (that is, there exists $l:[0,1]\to\real$ such that 
  $|\derpar\theta \widetilde C(\theta,v)|=l(v)$)
  \begin{proof}
    We just choose 
    \[ \varphi (\theta,v) = \frac{2\pi}{\len C}\int_0^\theta |\dot C(t,v)|~dt \]
  \end{proof}
\end{Proposition}

On the other hand, we may reparameterize to eliminate $\pi_T\derpar v C$
\begin{Proposition}\label{prop: black magic}
  For any $C^2$ regular homotopy $C$ such that all curves $\theta\mapsto C$
  are immersed
  there exists a $\varphi$ as in \eqref{eq: repar C} above 
  such that $\pi_{\widetilde T}\derpar v \widetilde C=0$.

  \begin{proof}
    Both  $\pi_T\derpar v   C$ and 
    $\overdot C\derpar v\varphi$ are parallel to $T$:
    so, if $\derpar v\varphi= 
      -{\langle\derpar v C, T\rangle}/
      {|\overdot C|}$, then 
      $\pi_{\widetilde T}\derpar v \widetilde C=0$
      
    The O.D.E.     
    \begin{equation}
      \begin{cases}
      \derpar v\varphi(\theta,v)= 
      -\frac{\langle\derpar v C(\tau,v), T(\tau,v)\rangle}
      {|\overdot C(\tau,v)|}
      \\
      \tau\defeq\varphi(\theta,v) \cr
      \varphi(\theta,0)=\theta, ~~~\theta\in S^1
    \end{cases}
    \label{eq:def_di_varphi_black}
  \end{equation}
    
    can be solved for $v\in[0,1]$, since 
    $$M(\tau,v)=-\frac{\langle\derpar v C(\tau,v), T(\tau,v)\rangle}
    {|\overdot C(\tau,v)|}$$
    is periodic in $\tau$ and continuous, and then is bounded: that is, 
    \[\max_{S^1\times[0,1]} M < \infinity\]

    Defining $\psi=\derpar \theta \varphi=\dot \varphi$  we compute
    \begin{eqnarray*}
      \der v\psi(\theta,v)= 
      \der \theta \der v \varphi=
      -\der \theta\frac{\langle\derpar v C(\tau,v) ,
      T(\tau,v)\rangle}{|\overdot C(\tau,v)|}
      =\\=
      -\psi(\theta,v)\der \tau\left(\frac{\langle\derpar v C(\tau,v) ,
      T(\tau,v)\rangle}{|\overdot C(\tau,v)|}\right)
    \end{eqnarray*}
    so $\psi$ solves
    \begin{equation}
      \label{eq: ode psi}
      \begin{cases}
        \derpar v\psi(\theta,v)= 
      -\psi(\theta,v)\der \tau\frac{\langle\derpar v C(\tau,v) ,
        T(\tau,v)\rangle}{|\overdot C(\tau,v)|}
        \\
        \psi(\theta,0)=1, ~~~\theta\in S^1
      \end{cases}
    \end{equation}
    and then $\derpar \theta \varphi>0$ at all times $v$.
  \end{proof}
\end{Proposition}
Note that $\varphi$ is not unique: we may change in 
\eqref{eq:def_di_varphi_black} and simply require that 
$\varphi(\cdot,0)$ be a diffeomorphism. \quad
The above result is stated in section 2.8 in \cite{Michor-Mumford}; there 
$\widetilde C$ is called a \emph{horizontal path}, since it
is the canonical parallel lifting of a path in $B_{i,f}$ to a path in $M$.

\smallskip

\begin{figure}[htbp]
  \centering
  \includegraphics[width=0.8\textwidth]{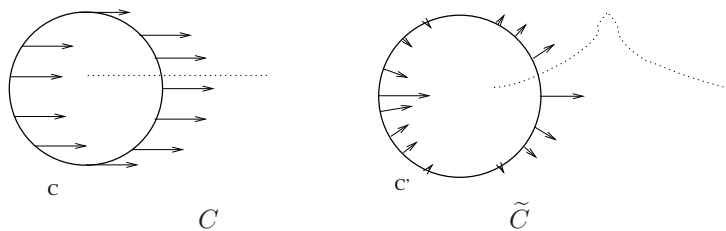}    

 \hfill {$C$} \hfill $\widetilde C$ \hfill~

  \caption[Reparameterization to $\pi_{\widetilde T}\derpar v \widetilde C=0$]
  {Reparameterization to $\pi_{\widetilde T}\derpar v \widetilde C=0$
    for the translation of the unit circle along the $x$-axis.
    The dotted line shows the trajectory of a point on the curve.}
  \label{fig:circle}
\end{figure}

For example, consider the unit circle translating with unit speed
along the $x$-axis, giving rise to the following homotopy $C(\theta,v)$:
\[
  C(\theta,v)=(v+\cos\theta,\sin\theta)
\]
While this is certainly the simplest way to parameterize the homotopy,
it does not yield a motion purely in the normal direction at each
point along each curve. However, the following reparameterization of
the homotopy yields exactly the same family of translating circles (and
therefore the same homotopy) but such that each point on each circle
along the homotopy flows exclusively along the normal to the corresponding
circle (i.e.\ in the radial direction).

\begin{eqnarray*}
  \widetilde C(\theta,v)&=&(\widetilde x,\widetilde y)
  \\
  \widetilde x(\theta,v)&=&v+\frac{(1-e^{2v})
                   +(1+e^{2v})\cos\theta}
                   {(1+e^{2v})+(1-e^{2v}) \cos\theta}
  \\
  \widetilde y(\theta,v)&=&\frac{2 e^v \sin\theta}
                   {(1+e^{2v})+(1-e^{2v}) \cos\theta}
\end{eqnarray*}
In figure \vref{fig:circle} we see a comparision between the
the trajectories of various points (fixed values of $\theta$)
along the translating circle in the original homotopy $C$ and
its reparameterization $\widetilde C$.

\smallskip

The above  result is suprising, kind of ``black magic'': it seems to
suggest that while modeling the motion of curves we can
neglect the tangential part of the motion $\pi_T\derpar v C$.
Unfortunately, however, this is not the case.
We begin by providing a simple example.
 
\begin{Example}
  The curve $C(\cdot,v)$ is translating as in figure \vref{fig: stretch},
  with unit speed:
  \[ C(\theta,v)\defeq c_0(\theta) + v e_1 \]
   
  \begin{figure}[htbp]
    \centering \includegraphics{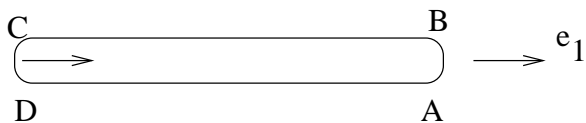}
    \caption{Stretching the parameterization}
    \label{fig: stretch}
  \end{figure}
  
  \noindent
  After the reparameterization \ref{prop: black magic},
  the points in the tracts $BC$ and $AD$
  are motionless, whereas the parameterization in $AB$ is stretching to
  produce new points for the curve, and in $CD$ it is absorbing
  points: then, $\derpar \theta \widetilde C\to\infinity$ if the
  curvature in $AB$ goes to infinity.
\end{Example}

We now consider the math in more detail.
\begin{Proposition}\label{no black magic}  
  Suppose that $n=2$ for simplicity: 
  define the curvature $\kappa$ as per definition \vref{def di k}.

  Suppose in particular that $|\overdot C|=1$ at all points: then
  \[ \derpar v |\overdot C|^2 =0= 2 \langle\derpar v\derpar \theta
  C , T\rangle\]
  
  By deriving  \eqref{eq: repar C} we obtain
    \begin{eqnarray*}
      \label{eq: d theta C   after black magic}
      \derpar \theta \widetilde C(\theta,v)=
      \derpar\theta C(\tau,v) \psi(\theta,v)
    \end{eqnarray*}
    where  $\psi=\derpar\theta\varphi$
    solves \eqref{eq: ode psi}, that we can rewrite
    (using $|\overdot C|=1$) as
    $$    \derpar v \psi=-\psi { \langle\derpar v C,N\rangle \kappa  }$$
    (where $C,\kappa ,N,T$ are evaluated at $(\varphi(\theta,v),v)$).

    Note that
    $|\langle\derpar v C, N\rangle|=
    |\pi_N\derpar v C|$.
    This implies that the parameterization of $\widetilde C$ will be highly affected
    at points where both $\kappa$ and $\pi_N\derpar v C$ are big.

    \label{rem: not both}
  So the above teaches us that there are two different approches to the
  reparameterization: we may either use it to control $\overdot C$ or
  $\pi_T\derpar v C$, but not both. 
\end{Proposition}

\subsection{Homotopy classes}

\newcommand\C{\mathbb C}
\subsubsection{Class $\C$}
Given an energy and two continuous curves $c_0$ and $c_1$,
we will try to find a homotopy that
minimizes this energy, searching  in the class $\C$
so defined

\begin{Definition}[class $\C$]
\label{class C}
  Let $I=S^1\times [0,1]$. Let $\C$ be the class of all homotopies
  $C:I\to\real^n$, continuous on $I$ and
  locally Lipschitz in $S^1\times(0,1)$, such that
  $c_0=C(\cdot, 0)$ and $c_1=C(\cdot, 1)$.
\end{Definition}
Such minimum will be a \emph{minimal geodesic} 
connecting $c_0$ and $c_1$ in the space of curves.
In this class $\C$ we can state 
\begin{Proposition}[Poincar\`e inequality]
  \label{poincare}
  there are two constants $a',a''>0$ such that $\forall C\in\C\cap H^{1,0}$,
  \begin{equation}\label{eq: poincare}
    \| C  \|_{H^{1,0}(I)} ^2\le a'+ a'' \int_I | \derpar v C | ^2   
  \end{equation}
  $a',a''>0$ depend on $c_0$ and $c_1$.
\end{Proposition}

\def\F{{\mathbb F}}
\subsubsection{Class $\F$ of prescribed-parameter curves}

In some of the following sections we will
change our point of view with respect to the above  assumption \ref{class C}
by fixing a measurable positive function $l:[0,1]\to \real^+$
and restricting our attention
to a family of homotopies such that $|\derpar \theta C(\theta,v)| =l(v) $.

\begin{Definition}[class $\F$]\label{class F}
  The class $\F$ contains all  of $C\in \C$ such that
  \[|\derpar \theta C(\theta,v)| =l(v)\] for all $\theta,v$.
\end{Definition}
In particular, $ \len(C)=2\pi l(v)$.

Note that, by \ref{prop: l is C}, if  $l$ is not H\"older continuous
then the class $\F$ will be too poor to be useful (for example,
it will not contain  smooth functions).
Unfortunately the class $\F$ is not closed with respect to the weak convergence
in $W^{1,p}$
\begin{Proposition}\label{prop: le l}
  Suppose that $l$ is bounded and 
  $|\derpar \theta C_h(\theta,v)| =l(v) $ for all $h,\theta,v$;
  assume that $1\le p<\infinity$ and 
  $\derpar \theta C_h \wto V$ weakly in $L^p(I\to\real^n)$, or
  $p=\infinity$ and   $\derpar \theta C_n \wto^* V$ weakly-* in
  $L^\infinity(I\to\real^n)$.

  Then   $|V(\theta,v)| \le l(v)$.
  \begin{proof}
    Since $|\derpar \theta C_h|\le \sup l$, we may always assume that
    $\derpar \theta C_h \wto^* V$ weakly-* in
    $L^\infinity(I\to\real^n)$; the result follows immediatly from
    theorem 1.1 in \cite{dacorogna:weak}
  \end{proof}
\end{Proposition}
We conclude that
\begin{Theorem}
  Let $1\le p \le \infinity$.
  The closure of \, $\F\cap W^{1,p}$ \, with respect to weak convergence $W^{1,p}$
  is contained in
  the class $\overbar \F$ of all $C\in W^{1,p}$ such that $C(\cdot,
  0)=c_0$ and $C(\cdot, 1)=c_1$ are given and $|\derpar \theta C| \le
  l(v) $. 
  \begin{proof}   
    Suppose $(C_h)\subset\F$ and $C_h\wto C$: we prove that $C\in\overbar \F$.
    $(C_h)$ is bounded in $W^{1,p}$: by
    Rellich-Kondrachov theorem (see thm.\ IV.16\ \cite{Brezis})
    $(C_h)$ is pre-compact in $L^p$: up to a subsequence, we may
    assume that $C_h\to C$ in  $L^p$ and that $C_h\to C$
    on almost all points: then $C(\cdot,0)=c_0$ and $C(\cdot,1)=c_1$.
    \ammark{On the other hand, suppose  $C\in\overbar \F$}
  \end{proof}
\end{Theorem}

\medskip 
\begin{Remark}
  
  We could just as well have defined a class of 
  \emph{homotopies with prescribed lengths}, where we fix
  a function $\hat l$ and include in the class all $C$ such that 
  $\len(C) =  \hat l(v)$ for all $v$. However, if the energy $E$ is geometric,
  then using the reparameterization \ref{prop: by arc par}, we can always
  replace any such $C$ with a $\widetilde C\in \F$, and
  $E(C)=E(\widetilde C)$. 
\end{Remark}

\subsubsection{Factoring out reparameterizations}

If we only consider curves $c$ such that $|\dot c|\equiv 1$, then
(as pointed in \cite{Michor-Mumford})
the group of reparameterizations is $S^1\ltimes \integer_2$, where
\begin{itemize}
\item the group  $S^1$ is associated to the change of the initial
  point  $c(0)$ for the  parameterization, that is, 
  the reparameterization $\tau\in  S^1$ changes  the curve $c$ to
  $c(\theta+\tau)$,
\item  and  $\integer_2$ means the operation of changing direction,
  that is, $c(\theta)$ becomes  $c(-\theta)$
\end{itemize}

We extend the above to homotopies satisfying 
$|\derpar\theta C(\theta,v)|=l(v)$;
the second reparameterization is not significant;  the first one
does not affect the normal velocity $\pi_N\derpar v C$;
so the reparameterization in the class $\F$
does not affect energies that depend only on $\pi_N\derpar v C$
(as per eq.~\eqref{eq: F = F pi h}).

\section{Analysis of $E^N(C)$}
\label{sec: ana}

In this section we will focus our attention on the energy 
\[ E^N(C) \defeq \int_{I}  |\pi_N \derpar v C|^2 |\dot C| \, d\theta \,dv  \]
which is associated with the geometric version of the $H^0$ metric.
We will derive a result of existence of minimima of $E^N$ in
a class of homotopies such that the curvature is  bounded.

\begin{Remark}[winding]\label{exa: E N not coerc}
  Unfortunately, any energy that is independent of parameterization
  cannot be coercive. It is then difficult to prove existence of geodesics
  by means of the standard procedure in the Calculus of Variations
  \ref{standard CV}.

  Indeed, suppose that $C$ is a homotopy, and define
  $$  C_k(\theta,v) = C(\theta+k2\pi v,v)$$  
  $\forall k\in\mathbb Z$; then
  \[ E^N(C_k) = E^N(C) \]

  (As a special case, consider two curves $c_0=c_1$, and
  $C_k(\theta,v)=c_0(\theta+k2\pi v)$; each $C_k$
  is a  minimal geodesic connecting $c_0,c_1$, since $E^N(C_k)=0$).

  This proves that $E^N$ is not coercive in $H^1(I)$, since
  $\int |\derpar v C_k|^2\to \infinity $ when $k\to\infinity$.
\end{Remark}

\subsection{Knowledge base}

Some of the follwing results apply
in the class $\F$ of homotopies with
prescribed arc parameter $|\derpar \theta C(\theta,v)|=l(v)$ while
others apply in the more general class $\C$.

\medskip

Consider the integral
\begin{equation}\label{eq:J(C)}
  J(C) \defeq \int_{I} |H|^2 |\pi_N \derpar v C|^2 |\dot C| \, d\theta\,dv 
\end{equation}
where $H$ is the curvature of $C$ along $\theta$.
We now deduce the following result from \cite{Michor-Mumford}
\begin{Proposition}\label{prop: l is C}
  Suppose that the homotopy $C$ is smooth and
  immersed: then by extending to $\real^n$ 
  the computation in sec.~3.3 of \cite{Michor-Mumford} we get
  \[ \der v \sqrt{\len(C)(v)} \le \frac 1 2 
  \sqrt {\int_{S^1} \langle H, \derpar v C\rangle^2 |\dot C|d\theta} \]
  and then, for any $0\le v'<v'' \le1$,
  \begin{eqnarray*}
    \sqrt{\len(C)(v'')}- \sqrt{\len(C)(v')} &\le&  \frac 1 2 
    \int_{v'}^{v''} \left(\int_{S^1}  \langle H, \derpar v C\rangle^2
      |\dot C|d\theta    \right)^{1/2} dv 
    \le \\ &\le&
    \frac 1 2     \sqrt{v''-v'} 
      \left( \int_{v'}^{v''}\int_{S^1}  \langle H, \derpar v C\rangle^2
      |\dot C|d\theta\,dv    \right)^{1/2} 
    \le \\ &\le&
    \frac{ \sqrt{ J(C) }}2 \sqrt{v''-v'}
  \end{eqnarray*}
  (by Cauchy-Schwarz and \eqref{eq:H,v<H,pi_nv})
  and this implies that $\sqrt {\len(C)}$ is H\"older continuous
  when  $J(C)$  is finite.
\end{Proposition}

\begin{Proposition}[l.s.c.\ and polyconvex function]\label{polyconvex}
  Let $W,V\in\real^n$.  Let $W\times V$ be the vector of all $n(n-1)/2$
  determinants of all 2 by 2 minors of the matrix having $W,V$ as
  columns.

  Consider a continuous function 
  $f:\real^{n\times 2}\to\real$ such that
  \[f(W,V)=g(W,V,W\times V)\]
  where   $g:\real^{n(n+3)/2}\to\real$ is convex: then $f$ is \emph{polyconvex}
  \footnote{The more general definition and the properties may be found 
    sec.\ 4.1 in \cite{buttazzo89:semic}, or
  in 2.5 and 5.4 in \cite{dacorogna:weak}}. 

  Let $p\ge 2$, suppose $f\ge 0$:
  by theorem 4.1.5 and remark 4.1.6 in \cite{buttazzo89:semic}
  then 
  \[ \int_I f(\derpar \theta C,\derpar v C) d\theta dv \]
  is $W^{1,p}$-weakly-lower-semi-continuous.\footnote{Sequentially 
    weakly-* if  $p=\infinity$}
  This means that,
  if  $ C_h\wto C$ weakly in $W^{1,p}$, that is
  \begin{equation}
    C_h\to C, ~~\derpar v C_h \wto  \derpar v C ~~
    \derpar \theta C_h \wto  \derpar \theta C\label{eq:vcrev}
  \end{equation}
  in $L^p$, 
  \footnote{We can write equivalently $C_h\to C$ or  $C_h\wto C$
    in \eqref{eq:vcrev}, thanks to
    Rellich-Kondrachov theorem (see thm.\ IV.16\ \cite{Brezis})}
  then     
  \[ \liminf_h \int_I f(\derpar \theta C_h,\derpar vC_h)
  \ge    \int_I f(\derpar \theta C,\derpar v C)  \] 
\end{Proposition}

\subsection{Compactness}
\label{sec: compactness}

We now list some simple lemmas.

\begin{Lemma} \label{prop: pi T con media}
  let $\widetilde C(\theta,v) =C(\theta+\varphi(v),v)$ be 
  a reparameterization;  then
  \begin{equation}
    \inf_\varphi \int |\pi_T \derpar v \widetilde C|^2=
    \int_0^1 \int_{S^1}\left ( \pi_T \derpar v  C (\theta,v)
      - \lint_{S^1} \pi_T \derpar v C (s,v)~ds
    \right)^2~d\theta~dv
    \label{eq: pi T con media}  
  \end{equation}
\end{Lemma}

\begin{Lemma}[Poincar\`e inequality in $S^1$]
  If $f:S^1\to\real^n$ is $C^1$, then  
  \begin{equation}
    \max |f|-\min |f| \le \frac 1 \pi  \int_{S^1} |d f|
  \end{equation}
  so that 
  \begin{equation}
    \sqrt{\int \left| f-\lint f \right|^2} \le \frac 1 {2\pi} \int_{S^1} |d f|
  \end{equation}
  For any $a\in\real, a\neq 0$, 
  \begin{equation}
    \max |f|-\min |f| \le   \frac 2 \pi  \int_{S^1} |d f+a|
    \label{eq: Poincare max min}
  \end{equation}
  so that 
  \begin{equation}
    \sqrt{\int \left| f-\lint f \right|^2} \le  \frac 1 \pi \int_{S^1} |d f+a|
    \label{eq: Poincare L^1}
  \end{equation}
\end{Lemma}

\begin{Lemma}
  \label{derivatives}
  Suppose $C\in C^2$.
  If we derive $\derpar v(|\overdot C|^2)$ we get
  $$\derpar v(|\overdot C|^2)= 2\langle T ,\derpar v\derpar\theta C
  \rangle |\overdot C|$$
  so,%
  \footnote{  When $|\overdot C|=0$, then $T=0$, by our definition, 
    so the equation \eqref{eq:der_di_arco} holds as well
    in the distributional sense}
  when $|\overdot C|\neq 0$,
  
  \begin{equation}
    \label{eq:der_di_arco}
    \derpar v (|\overdot C|)= \langle T ,\derpar v\derpar\theta C  \rangle    
  \end{equation}  

  If the curves are in $C^2$ and
  we derive $\derpar\theta (\pi_T \derpar v C)$, we get  
  \begin{equation}
    \label{eq:der_di_pi_N_d_v}
    \derpar\theta (\pi_T \derpar v C)= 
    \derpar\theta\langle T ,\derpar v C \rangle= 
    \langle H ,\derpar v C\rangle|\overdot C|
    + \langle T,\derpar \theta\derpar v C  \rangle= 
    \langle H ,\derpar v C\rangle|\overdot C|
    +   \derpar v (|\overdot C|)
  \end{equation}
  (where $\langle H ,\derpar v C\rangle=  \kappa \pi_N\derpar v C$
  for curves in the plane).
\end{Lemma}

\begin{Proposition}\label{prop: bound H1}
  Let $M>0$ be a constant.
  Suppose that a  $C^2$ homotopy  \hbox{$C:I\to\real^n$} satisfies  
  $|\dot C(\theta,v)|=l(v)$ and
  \begin{equation}
    \int_0^1 \left(  \int_{S^1}
      \langle H, \derpar v C \rangle |\dot C|~d\theta\right)^2dv \le M
    \label{eq: bound on H Cv}
  \end{equation}
  where $H$ is the curvature of $C$ along $\theta$.
  Then there exists a suitable reparameterizations   
  \begin{equation}
    \widetilde C(\theta,v)= C(\theta+ \varphi(v),v)\label{eq:unwinding}    
  \end{equation}
  such that
  \[ \int |\pi_T \derpar v \widetilde C|^2 \le 2M \]

  \begin{proof}    
    Suppose that $l\in C^1$ (the general proof being obtained by an
    approximation argument).
    We summarize  derivations in \ref{derivatives},
    \[\derpar\theta (\pi_T \derpar v C)=  
    \langle H,\derpar v C \rangle |\dot C|+
    \derpar v l\]
    Applying eqns   \eqref{eq: pi T con media}   and \eqref{eq: Poincare L^1}
    (with $a=-\derpar v l$) we obtain
    \begin{eqnarray*}
      \inf_\varphi \int |\pi_T \derpar v \widetilde C|^2=
      \int_0^1 \int_{S^1}\left | \pi_T \derpar v  C (\theta,v)
        - \lint_{S^1} \pi_T \derpar v C (\tilde \theta,v)~d\tilde\theta
      \right|^2~d\theta~dv
      \le \\ \le
      \int_0^1 \left(  \int_{S^1}
        |\derpar\theta (\pi_T \derpar v C)+a| ~d\theta  \right)^2
       \le
      \int_0^1 \left(  \int_{S^1}
        |\langle H, \derpar v C\rangle |\dot C| 
         +\derpar v l-\derpar v l | ~d\theta\right)^2dv
    \end{eqnarray*}

  \end{proof}
\end{Proposition}
Note that the reparametrization \eqref{eq:unwinding} may be viewed as an 
\emph{unwinding} when confronted with \ref{exa: E N not coerc}.

\begin{Remark}\label{rem: bound on J compactness}1
  Note that, if $|\dot C(\theta,v)|=l(v)$ then
  \[   \int_0^1 \left(  \int_{S^1}
    \langle H, \derpar v C \rangle |\dot C|~d\theta\right)^2 dv \le   
  \int_I \langle H, \derpar v C \rangle^2 l(v)^2\le   (\sup l) J(C)  \]
  So a bound on $J(C)$ should provide compactness.
\end{Remark}

\ammark{The proposition \ref{no black magic} shows that there is no
better reparameterization than the one chosen in the theorem.
}

\subsection{Lower semicontinuity}
\label{sec: lsc}

\begingroup
\def\a{\alpha}
\def\b{\beta}
\def\t{\theta}

Let $\a>0, \b>0$, $V,W\in\real^n$,  define  
\[ e(W,V) = |\pi_{W^\perp}V |^\a |W|^\b \]
and define
\[E_{\a,\b}(C) \defeq   \int_I e(\derpar\t C,\derpar v C) \] 
Note that $E^N(C)$ is obtained by choosing $\a=2,\b=1$.
In general if $\beta=1$ then $E_{\a,\b}(C)$ is a geometric energy
(see \ref{def: geo ene}).

Let $W\times V$ be the vector of all $n(n-1)/2$
determinants of all 2 by 2 minors of the matrix having $W,V$ as
columns.
The identity 
\[ |\pi_{W^\perp}V |^2 |W|^2 =
  |V|^2|W|^2 -\langle  V,W\rangle^2=|V\times W|^2 \]
is easily checked 
\footnote{In $\real^3$ we have
  $\langle V,W\rangle=|V||W|\cos\alpha$ and $|V\times W|=|V||W|\sin\alpha$,
  where $\alpha$ is the angle between the two vectors}.
Let   
\[ f(W,V)=|V\times W| ~~. \]
Note that $f$ is a polyconvex function (see \ref{polyconvex})
and that
\begin{equation}
  e(W,V) = |\pi_{W^\perp}V |^\a |W|^\b = |W\times V|^\a |W|^{\b-\a} 
  \label{eq:|VxW|}
\end{equation}
and
\[ E_{\a,\b}(C)= \int_I f(\derpar \theta C,\derpar v C)^{\a}
|\derpar\theta C|^{\b-\a} \]

\smallskip

We can provide this lower semicontinuity result in the class $\F$.

\begin{Proposition}\label{lsc}
  Let $p\ge 2$.
    Suppose $\a> \b>0$.
    Fix a continous function $l:[0,1]\to\real^+$, $l\ge 0$, and use it to
    build the class $\F$.  
    Then $E_{\a,\b}$ is $W^{1,p}$-weakly-lower-semi-continuous in the class
    $\F$. 

    More precisely: let
    \[ C\in \F, ~~~~ (C_h)_h\subset \F  \]
    if  $ C_h\wto C$ weakly in $W^{1,p}$, that is,
    \[ C_h\to C, ~~\derpar v C_h \wto  \derpar v C ~~
    \derpar \theta C_h \wto  \derpar \theta C   \]
    in $L^p$,  and 
    \[ l(v) =  |\dot C(\theta,v)| = |\dot C_h(\theta,v)|   \]
    then     
    \[ \liminf_h E_{\a,\b} (C_h) \ge E_{\a,\b}(C) \]
\end{Proposition}
\begin{proof}
  We prove %
  the theorem in steps:
  \begin{itemize}
  \item 
      Let $\lambda>0$.
      Suppose  %
      $|\dot C|\equiv |\dot C_h|\equiv \lambda$; then 
      \[ E_{\a,\b}(C)= \lambda^{\b-\a} \int f(\derpar \theta C,\derpar v C)^\a \]
      is l.s.c.\ (by \ref{polyconvex}).

    \item
      for any homotopy $C$ and any continuous $g:[0,1]\to\real^+$,
      define $e_g(C)(v):[0,1]\to \real^+$,
      \[ e_g(C)\defeq \oldcases{ 
         {g(v)^{\b-\a}}
        \int_{ S^1}     f(  \derpar \theta C  ,\derpar v C)^\a d\theta
        & if $g(v)>0$ \cr
        0& if $g(v)=0$} \]
      
      If $C\in \F$ (that is   $l(v)=|\dot C(\t,v)|$)  then
      \[ E_{\a,\b}(C)    =  \int_0^1  e_l(C) dv\]

    \item 
      Consider a piecewise function  $g\ge 0$ defined 
      \begin{equation}
        g=\sum_{i=1}^m  g_i \chi_{[a_i,a_{i+1})}\label{eq:def-di-pw-g}        
      \end{equation}

      and let 
      \[ \hat E(C) \defeq   \int_0^1  e_{g}(C) dv =
      \sum _{i=1}^m       \int_{a_i}^{a_{i+1}}      e_{g_i}(C) d v  \]

      then we apply the previous  reasoning to all addends
      and conclude that 
      \[      \liminf_h \hat E(C_h) \ge  \hat E(C)\]

    \item 
      Suppose $l$ is continuous and $l\ge 0$.
      Let $\tau$ be the class of piecewise functions  $g\ge 0$ defined 
      as in \eqref{eq:def-di-pw-g},
      such that on any interval $[a_i,a_{i+1})$, either 
      $ g_i=0$ or 
      \[ g_i\ge \sup_{[a_i,a_{i+1})} l ~~.\]

      Then for any such $g$ and $C\in \F$,
      \[ \int_0^1  e_{g}(C) dv \le E_{\a,\b}(C) ~~.\]
    \item 
      choose $g$ in the class $\tau$;      then
      \[ \liminf_h E_{\a,\b}(C_h) \ge\liminf_h  \int_0^1  e_{g}(C_h) dv
      \ge  \int_0^1  e_{g}(C) dv  \]

    \item 
      Fix $C$: it is possible to find a sequence $(g_j)\subset \tau$ 
      such that  
      \[ e_{g_j}(C)(v) \to_j e_l(C)(v) \]
      for almost all points  $v$, monotonically increasing: 
      indeed, let
      \[A_{j,i} = [i2^{-j},(i+1)2^{-j}), ~~~
      I_{j,i} = \inf_{v\in A_{j,i}} l(v),~~~
      S_{j,i} = \sup_{v\in A_{j,i}} l(v)\]
      and
      \[ g_j(v) = \oldcases {
        0 & if $v\in A_{j,i}$ and $ I_{j,i}=0$        \cr
        S_{j,i}& if $v\in A_{j,i}$ and $ I_{j,i}>0$       }\]

      Then 
      \[  E_{\a,\b}(C) = \sup_j  \int_0^1  e_{g_j}(C) dv      \]
    \end{itemize}
  \end{proof}

We would like to prove this more  general statement
\begin{Conjecture}\label{conj:sci}
 choose Lipschitz homotopies
  \[ C:I\to\real^n, ~~ C_h:I\to\real^n \]
  define
  \[ l(v) \defeq \len C, ~~ l_h(v) \defeq \len C_h   \]
  
  If $l_h\to l$    uniformly and
  $ C_h\wto C$ weakly in $W^{1,p}$,  then 
  
  \[ \liminf_h E^N (C_h) \ge E^N(C) \]
\end{Conjecture}

Whereas we cannot generalize the theorem
further: this is due to example \ref{exa: E not lsc}.

\endgroup

\subsubsection{Example}
\label{sec: tesselation}
We want to show that $E^N$ is not l.s.c.\ 
if we do not control the length.

\begingroup
\label{exa: E not lsc}
\def\a{\alpha}
\def\b{\beta}
\def\t{\theta}

Let $\a>\b>0$, define  
\[ e(W,V) = |\pi_{W^\perp}V |^\a |W|^\b = |W\times V|^\a |W|^{\b-\a} \]
and define
\[E_{\a,\b}(C) \defeq   \int_I e(\derpar\t C,\derpar v C) \] 

We  will actually show that $E_{\a,\b}$ is not l.s.c., and that
\begin{Proposition}\label{prop: relaxed 0}
  \[\Gamma E_{\a,\b}(C)=0\]
  where the relaxation is computed with respect to weak-* $L^\infinity$
  convergence of the derivatives.
\end{Proposition}

\begin{enumerate}
\item 
 For simplicity, we temporarily drop the
  requirement that the curves $C(\cdot,v)$ be closed. 
  We then redefine $I=[0,1]\times[0,1]$.

\item 
   Let  $\widetilde C:I\to I$ be a Lipshitz map such that
   \begin{equation}
     \widetilde C(\t,0)= (\t,0), ~~ \widetilde C(\t,1)= (\t,1), ~~
     \widetilde C(0,v) = (0,v), ~~\widetilde C(1,v) = (1,v) ~~.\label{eq:d3409fh}
   \end{equation}

 \item \label{item:tessel}
  Let $h\ge 1$ be integers.
 
 We rescale $\tilde C$ and glue many copies of it to build $C_h$, as
 follows
 \[ C_h(\t,v) \defeq 
 \frac 1 h \tilde C\Big( (h\t)\mod (1), (h v)\mod(1)\Big) 
 + b_h(\t,v) \]
 where $b_h:I\to\real^2$ is the piecewise continuous function
 \[ b_h(\t,v) \defeq \left ( 
   \frac{1} h \lfloor h\t \rfloor ,
   \frac{1} h \lfloor h v \rfloor 
 \right ) \]
 (in particular, $C_1=\tilde C$).
 We represent this process in figure \vref{fig:texel}.
 \begin{figure}[htbp]
  \centering
  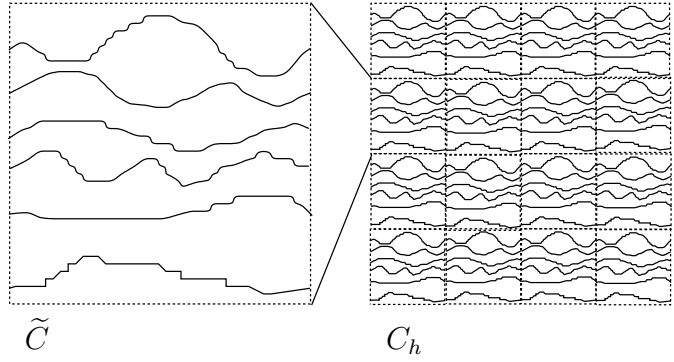
  \caption{Tesselation of homotopy $\tilde C$    to form $C_h$}
  \label{fig:texel}
\end{figure}

 Then
 \begin{eqnarray*}
   \derpar \t C_h(\t,v) = & \derpar \t \tilde C\Big( (h\t)\mod (1), (h v)\mod(1)\Big) \\
   \derpar  v C_h(\t,v) = & \derpar  v \tilde C\Big( (h\t)\mod (1), (h v)\mod(1)\Big)
 \end{eqnarray*} 

 We may think of $C_h$ as 
 homotopies that connect the same two curves, namely,
 \[ C_h(\t,0)= (\t,0)= \tilde C(\t,0),~~~C_h(\t,1)= (\t,1)= \tilde C(\t,1)~~,\]
 while extrema points move in a controlled way, namely
 \[ C_h(0,v)= (0,v)= \tilde C(0,v),~~~C_h(1,v)= (1,v)=\tilde C(1,v) ~~.\]

\item 
 Let $C(\t,v)\defeq (\t,v)$ be the identity.

 The sequence $C_h$ has the following properties
 \begin{enumerate}
 \item $C_h\to C$ in $L^\infinity$, and more precisely
   \[ \sup_{I} |C_h-C| \le \frac 2 h \]

 \item $\derpar \t C_h$ and $\derpar v C_h $
   are bounded in $L^\infinity$, and more precisely,
   \[ \sup_I |\derpar \t C_h| \le  |\derpar \t \tilde C|,~~~
   \sup_I |\derpar v C_h|\le  \sup_I |\derpar v \tilde C| \] 
   and then all $C_h$ are equi Lipschitz;

 \item \[ \derpar \t C_h \wto e_1=(1,0)=\derpar \t C \] 
   and \[ \derpar v C_h \wto e_2=(0,1)=\derpar v  C \]
   weakly* in $L^\infinity(I)$.~%
\footnote{Proof by lemma 1.2 in \cite{dacorogna:weak}}
   
 \item let $\tilde l(v)\defeq \len(\tilde C)$, $ a\defeq \int \tilde l$ and
   \[ l_h(v) \defeq \len(C_h)(v) = \tilde l\big( (hv)\mod(1)\big)\]
   then $l_h\wto a$   weakly* in $L^\infinity([0,1])$.$~^{\thefootnote}$

 \item\label{item: J esplode}
   Suppose $\tilde C$ is piecewise smooth, so that the curvature
   $H_h$ of $C_h$ can be defined almost everywhere.
   By \eqref{eq:d3409fh}, $\tilde l(0)=\tilde l(1)=1$.
   If $\tilde l(v)>1$ at some points, then the sequence
   $l_h(v)$ is not equicontinuous: then, by 
   proposition~\ref{prop: l is C}, the sequence of integrals
   \[ \int_{0}^{1}\int_{S^1}  \langle H_h, \derpar v C_h\rangle^2
   |\dot C_h|d\theta\,dv \]
   is unbounded in $h$.\footnote{a fortiori,
     $J(C_h)$ is unbounded --- $J(C)$ was defined in \eqref{eq:J(C)}}
 \item
   $E_{\a,\b}(C_h)$ is constant in $h$:
   indeed
   \begin{eqnarray*}
     \int_I e(\derpar\t C_h,\derpar v C_h)&=&
     h^2 \int_0^{1/h}\int_0^{1/h} 
     e\big(\derpar\t C_h,\derpar v C_h\big) d\t\,dv=\\
     &=& h^2\int_0^{1/h}\int_0^{1/h} 
     e\big(\derpar\t \tilde C(h\t,hv),\derpar v\tilde C(h\t,hv)\big) d\t\,dv=\\
     &=&\int_I e\big(\derpar\t \tilde C(\t,v),\derpar v \tilde
     C(\t,v)\big)   d\t\,dv
   \end{eqnarray*}
   that is 
   \begin{equation}
     \label{eq: E a b const}
     E_{\a,\b}(C_h)= E_{\a,\b}(\tilde C)
   \end{equation}
\end{enumerate}

\item
  Let $j\ge 1$ be a fixed integer.
  Let  $\widetilde C:I\to\real^2$ be defined by
  \begin{equation}
    \gamma(v)=\oldcases{
      v  & if  \(v\in[0,1/2]\) \cr
      (1-v) & if  \(v\in[1/2,1]\)
    }
  \end{equation}
  and
  \begin{equation}
    \tilde C(\t,v) =\big( \t,  v+\gamma(v)\sin(2\pi j\t) \big)
  \end{equation}
  (note that $\tilde C$ is a graph).

  We have this further properties
 \begin{itemize}
 \item 
   We know that  $l_h\wto a$ with $a=\int \len \tilde C$;
   but this limit $a$ is strictly bigger than $1$, whereas 
   $\len   C\equiv 1$
 \item
   For $\t\in[0,1/h]$ and $v\in[0,1/h]$
   \begin{eqnarray*}
     \derpar  v C_h(\t,v) =&  \derpar  v \tilde C(h\t,hv)
     &=(0,~ 1+\gamma'(hv)\,\sin(2\pi jh\t))\\
     \derpar \t C_h(\t,v) =&  \derpar \t \tilde C(h\t,hv)
     &=(1,~ 2\pi j  \gamma(hv)   \cos(2\pi jh \t) )
   \end{eqnarray*} 
   
   Then
   \[ \sup_I |\derpar \t C_h| \le 1+2\pi j ,~~
   \sup_I |\derpar v C_h|\le 2   \] 
   
 \item 
   We compute
   \begin{eqnarray*}
     E_{\a,\b}(\tilde C) =
     \int_I e\big(\derpar\t \tilde C(\t,v),\derpar v \tilde
     C(\t,v)\big)   d\t\,dv = 
     \\ =
     \int_0^{1} \int_0^1 |1+\gamma'(v)\sin(j\t)|^{\a}
     |1+\gamma(v)^2 j^2\cos(j\t)^2   |^{(\b-\a)/2} \,d\t\,dv
     \\ =
     2\int_0^{1/2} \int_0^1 |1+\sin(j\t)|^{\a}
     |1+v^2 j^2\cos(j\t)^2   |^{(\b-\a)/2} \,d\t\,dv
   \end{eqnarray*}
   then
   \begin{equation}     \label{eq: E a b j}
     \lim_{j\to\infinity} E_{\a,\b}(\tilde C) = 0
   \end{equation}
 \end{itemize}
\item 
  Combining \eqref{eq: E a b const} and \eqref{eq: E a b j}  
  we prove that  $E_{\a,\b}$ is not  l.s.c.\ Indeed for
  $j$ large
  \[ \lim_h E_{\a,\b}(C_h)=  E_{\a,\b}(\tilde C)
  <    E_{\a,\b} (C)=1 \]
  whereas $C_h\wto C$. 
\item to prove \ref{prop: relaxed 0}, consider any homotopy $C$;
  this may be approximated by a piecewise linear  homotopy,
  which in turn may be approximated by many replicas of the above construction.
\end{enumerate}

\endgroup

\subsection{Existence of minimal geodesics}

\begin{Theorem}\label{thm:existence}
  Let $M>0$. 

  Let $\mathcal A$ be the set of admissible curves $c:S^1\to\real^n$, such that
  \begin{itemize}
  \item $c:S^1\to\real^n$ is Lipschitz, and $c$ admits curvature  $H$
    in the measure sense (see in \ref{def di H}), and moreover
  \item the total mass $|H|(S^1)$ of the curvature $H$ of $c$
    is bounded uniformly
    \begin{equation}\label{eq: bound mass H}
       |H|(S^1) \le M ~~~.
    \end{equation}
  \end{itemize}

  Let $c_0,c_1$ be curves in  $\mathcal A$.

  Fix a bounded continuous function $l:[0,1]\to\real^+$, with  $\inf l>0$.

  Let $\mathcal B$ be  the class
  of homotopies $C:I\to \real^n$ such that
  \begin{itemize}
  \item $C\in H^1(I\to\real^n)$
  \item 
    any given curve $\theta\mapsto C(\theta,v)$ is in $\mathcal A$ 
  \item \label{cond compatt}
    the curvature can be extended to the homothopy (see \ref{extend H}),
    $\derpar v C$ is continuous, and
    \[ \int_0^1 \left(  \int_{S^1}
    \langle H, \derpar v C \rangle |\dot C|~d\theta\right)^2 dv   
    \le M   \tag{\eqref{eq: bound on H Cv}}    \]
  \item $l(v)=\len C(v)=\int_{S^1}|\dot C|~d\theta $ 
  \item    $C(\theta,0)=c_0(\theta)$ and $C(\theta,1)=c_1$
  \end{itemize}

  If $\mathcal B$ is non empty, then the functional $E^N$ admits a minimizing
  homotopy $C^*$;
  this minimum $C^*$ satisfies all the above requirements, but possibly for
  condition \eqref{eq: bound on H Cv}.
\end{Theorem}

\begin{proof}
  The proof is divided in two important (and independent) steps
  \begin{itemize}
  \item 
    Let $C_h$ be a sequence such that
    \[ \lim E^N(C_h) = \inf_{C\in \mathcal B} E^N(C)\]

    Up to reparameterization, assume
     \[ |\dot C_h(\theta,v)|=l(v) ~~. \]

    By this bound, and the compactness result \ref{prop: bound H1},
    we reparameterize any term $C_h$ to $\widetilde C_h$ by
    \[ \widetilde C_h(\theta,v) =C_h(\theta+\varphi_h(v),v)\]
    so that 
    \[ \int |\pi_T \derpar v \widetilde C_h|^2 \le 2M ~~~;\]
    moreover 
    \begin{equation}
      \int |\pi_N \derpar v \widetilde C_h|^2 \le
      (\max \frac 1{l})\int |\pi_N \derpar v \widetilde C_h|^2 |\dot C|\le
      (\max \frac 1{l})  (1+ \inf E^N) 
      \label{eq:acci} 
    \end{equation}
    (definitively in $h$) and then 
    \[ \int | \derpar v \widetilde C_h|^2 \le 2M+ 
    (\max \frac 1{l})   (1+\inf E^N) \]
    whereas
    \[ \int_I |\dot C_h|^2 = \int_0^1 l(v) ^2 dv~~~:\]
    then (by the Banach-Alaoglu-Bourbaki theorem
    \footnote{See thm III.15 III.25 and cor III.26 in \cite{Brezis}})
    up to a subsequence,  $\widetilde C_{h}$ converges
    weakly in $H^1$ to a homotopy $C^*$. 

  \item 
    We want to prove that 
    $|\dot C^*(\theta,v)|=l(v)$ for almost all $\theta,v$.
    We know that $|\dot C^*(\theta,v)|\le l(v)$,
    by \ref{prop: le l}. Suppose on the opposite that 
    $|\dot C^*(\theta,v)|<l(v)$: then there exists $\e>0$ and
    an measurable subset $A\subset I$ with positive measure, such that    
    $|\dot C^*(\theta,v)|<l(v)-\e$ for $(\theta,v)\in A$. 
    Let $A_v=\{ \theta : (\theta,v)\in A\}$ be the slice of $A$:
    then, by Fubini-Tonelli, there is a $v$ such that the measure
    of $A_v$ is positive. We fix that $v$. Suppose that $H_h$
    is the curvature of $C_h(\cdot,v)$: then 
    \[ H_h = l(v) \derpar{\theta\theta} C_h \]
    in the sense of measures.
    We know that $H$ has bounded mass: so $\derpar{\theta\theta} C_h$ has:
    by Theorem 3.23 in \cite{AFP:BV}, $\derpar{\theta} C_h(\cdot,v)$
    is compact in $L^1(S^1)$, so,   (up to a
    subsequence), we would have that $\dot C_{h}(\cdot,v) \to \dot
    C^*(\cdot)$ \emph{strongly} in $L^1(S^1)$:
    then $|\dot C^*(\theta,v)|=l(v)$ for that particular choice of $v$,
    achieving contradiction.
  \item 
    By \ref{lsc},
    $\liminf_h E^N (C_{k_h}) \ge E^N(C^*)$: then $C^*$ is the minimum.
\end{itemize}
\end{proof}

\begin{Remark}\label{rem: extend thm}
  If we wish to extend the above theorem, we face some obstacles.
  \begin{itemize}
  \item If we do not enforce some  bounds on curvature
    (as \eqref{eq: bound mass H} and
    \eqref{eq: bound on H Cv}), then
    the example in \S\ref{pulley} shows that we cannot achieve
    compactness of a minimizing subsequence 

  \item 
  If we wish to remove the hypothesis ``$\inf l > 0$'',%
\footnote{note that the bound \eqref{eq: bound mass H}
  does not imply that $\inf l > 0$}
  we are faced with the following problem: if $l(v)=0$, then 
  the curve $C(\cdot,v)$ collapses to a point; consequently
  if $l(v)=0$ on an interval $[a,b]$,
  then the homotopy collapses to path on that
  interval; moreover, the length and the energy of $C$
  restricted to $v\in[a,b]$ is necessarily $0$, so that $E(C)$ provides
  no bound on the behaviour of $C$: we again lose 
  compactness of a minimizing subsequence
  (indeed, the inequality 
  \eqref{eq:acci} needs that  $\inf_v\len(C)(v)>0$, to be able to
  control $\int |\derpar v C|^2$).
  \end{itemize}
\end{Remark}

\subsubsection{Michor-Mumford}
\label{sec: discuss MM}

Let $d(c_0,c_1)$  be the geodesic distance induced by 
the metric  $G^A$  defined in \cite{Michor-Mumford}
(see eq.\eqref{eq:G MM} here).
By the results in  \cite{Michor-Mumford}, we know that (in 
general) this distance is non degenerate:
\begin{Proposition}\label{rem: d>0 in MM}
  Consider an homotopy $C$ connecting two curves 
  $c_0=C(\cdot,0),c_1=C(\cdot,1)$,
  and its energy $E^A(C)\defeq E^N(C)+AJ(C)$
  (see eq.(\ref{eq: E da MM}) here)

  By H\"older inequality,
  \begin{eqnarray*}
    \left(\int_I |\dot C|\right)^{1/2}
    \left(\int_I |\pi_N\derpar v C|^2|\dot C|\right)^{1/2}\ge
    \int_I |\pi_N\derpar v C| |\dot C|\ge
    \int_I |\derpar v C\times \dot C|
  \end{eqnarray*}
  We then obtain the
  \emph{area swept out bound} 
  \footnote{in a form slightly better than the one 
    in \S3.4 in \cite{Michor-Mumford}}
  \[ \left(\int_I |\derpar v C\times \dot C| \right)^{2}
  \le \frac{E^N(C)}{\int_0^1\len(C)} \]
  Indeed the leftmost term is the area swept by the
  homotopy.

  By the proposition \ref{prop: l is C}, we know that, if $J(C)$ is
  bounded, then  $\len(C)$ is continuous and is bounded.
  So, if $c_0\neq c_1$, and there is no zero-area homotopy
  \footnote{note that there are  curves $c_0\neq c_1$
    in the completion of the
    space $B_i$ (the completion is described in \S3.11 in
    \cite{Michor-Mumford}) that can be connected by a   zero-area
    homotopy}
  connecting $c_0,c_1$, then $d(c_1,c_2)>0$. 
  See also prop.~\ref{prop: conforme non deg} here.
\end{Proposition}

We cannot instead currently prove a theorem of existence of
minimal geodesics for the energy $E^A(C)$;
in this section we have though derived some insight;
so we discuss the conjecture 

\begin{Conjecture}%
  \label{conj:MM}
  Fix two curves $c_0$ and $c_1$.
  The energy $E^A(C)$ admits a minimum in the class $\C$ of
  homotopies connecting $c_0$ and $c_1$.
\end{Conjecture}

May we improve the proof of theorem~\ref{thm:existence}
to prove this conjecture? We discuss what is 
ok and what is wrong.
\begin{itemize}
\item We may want to use $J(C)$ to drop the requirement~\eqref{eq: bound mass H}
  that is, in turn, used to have l.s.c.\ of the functional $E^N$;
  so we may think of proving this lemma
  \begin{center}
    \emph{``Suppose that we are given a sequence of
      smooth homotopies $C_h$, 
      and $J(C_h)\le M$, and $C_h\wto C$ weakly in $W^{1,p}$: 
      then $\len(C_h)\to \len(C)$''}
  \end{center}
  but this is wrong, as seen by this example
  \begin{Example}\def\t{u}
    Let $C_h:[0,1]\times[0,1]\to \real^2$ defined as
    \[ C_h(\t,v) = \left (\t,~ \frac 1 h \sin(2\pi h\t) \right) \]
    and
    \[ C(\t,v) = (\t,0)~~.\]

    These homotopies 
    do not depend on $v$: then 
    $J(C_h) = 0$. On the other hand, $C_h\wto C$ but $\len(C_h)$ is
    constant and bigger than $1=\len(C)$.
  \end{Example}
  
\item We need way to be sure that curves in the homotopy do
  not collapse to points (as discussed in \ref{rem: extend thm});
  so we may think of proving this lemma
  \begin{center}
    \emph{``Suppose that the homotopy $C$ admits curvature,
      and $J(C)<\infinity$: then $\inf_v\len(C)(v)>0$''}
  \end{center}
  but this is wrong, as seen by this example
  \begin{Example}\def\t{\theta}
    Let \[c_1(\t)=(\sin(\t),\cos(\t) )\]
    be the circle in $\real^2$, and build the homotopy
    \[ C(\t,v)=v^4 c_1(\t)\]
    Then
    \[ J(C) = \int_0^1 \frac 1{v^8} (4v^3)^2  2\pi v^4 dv =
    32\pi \int_0^1 v^2dv=
    32\pi/3    \]
  \end{Example}

\item We need a semicontinuity result on $J(C)$.
  Indeed we cannot hope in any cancellation effect
  in the sum $E^N(C)+J(C)$ because the two energies scale
  in different ways:
  \footnote{and this suggests that the  energy $E^A(C)$ 
    should not satisfy the rescaling property  \vref{pro: rescaling}}
  \begin{Remark}[rescaling]
    Let $\e>0$ and 
    $\tilde C= \e C$ then $E^N(\tilde C)=\e^3 E^N(C)$
    but  $J(\tilde C)=\e J(C)$
  \end{Remark}

\item On the bright side, we do not have a counterexample
  to show if $J(C)$ is not l.s.c.; actually,
  remark (\ref{item: J esplode}) (on page \pageref{item: J esplode})
  suggests that if  $C_h\to C$ and the homotopies have a common border 
  condition (such as \eqref{eq:d3409fh}),
  and $J(C_h)$ is  bounded, then $\liminf E^N(C_h)\ge  E^N(C)$.

\item 
  Moreover, by \ref{rem: d>0 in MM} we know that 
  the induced distance is in general not degenerate.

\item As pointed in \ref{rem: bound on J compactness},
  the term $J(C)$ in the metric provides compactness in $H^1(I)$

\item Moreover the  example in \S\ref{pulley} shows 
  that we do need to control the curvature of curves to be able
  to prove existence of minimal geodesics: this justifies the term $J(C)$
  in the Michor-Mumford energy $E^A$ (as
  well as the bound \eqref{eq: bound on H Cv} in theorem \ref{thm:existence}).
\end{itemize}

\subsection{Space of Curves}
\label{sec: space of bounded curvature}
By using the previous theorem~\ref{thm:existence}, we immediatly obtain 
a metric on a Space of Shapes.

Fix $M>0$.
Let $\mathcal S$ be the space of unit length
\footnote{if the conjecture \ref{conj:sci} holds true, then
 the ``unit length'' constraint may be dropped; 
 note that the formula of the metric is ``geometric'' (as defined in 
 \ref{def: geo ene}); but the minimal geodesics and the
 distance are not invariant
 with respect to rescaling, due to the bound $|\kappa| \le M$.
} closed immersed $C^2$ curves such
that for any $c\in \mathcal S$, the curvature $\kappa$ of $c$ is
bounded by $M$, as 
\[ |\kappa| \le M ~~.\]

We may think of  $\mathcal S$ as a ``submanifold with border'' in the
manifold $M$ of all closed unit length immersed $C^2$ curves.

Then we can use 
the Riemannian metric
\[ \langle h,k\rangle =
\int_{S_1} \langle \pi_Nh(\theta),\pi_N k(\theta)\rangle ~d\theta \]
to define a positive geodesic distance  in $\mathcal S$:
by the theorem in the previous section, this distance admits minimal
geodesics:
\begin{proof}
  indeed, we fix $\len(C)\equiv 1$;  we may write 
  \[ J(C) \le K^2 E^N(C) \]
  then we may use the remark  \ref{rem: bound on J compactness}
  to obtain compactness of minimizing subsequences.
\end{proof}

Unfortunately, since  $\mathcal S$ has a border (given by the
constraint $|\kappa| \le M$) then the minimal geodesic will
not, in general, satisfy the Euler-Lagrange ODE defined by $E^N$.

\subsection{The pulley}
\label{pulley}

We show that there exists a sequence of Lipschitz
functions  $C_h:I\to\real^2$ such that $|\pi_N \derpar v C_h|\le 1$,
$|\derpar\theta C_h|=1$,
but $\int_I |\derpar v C_h|\to \infinity$.

\begin{figure}[htbp]
  \centering
  \includegraphics{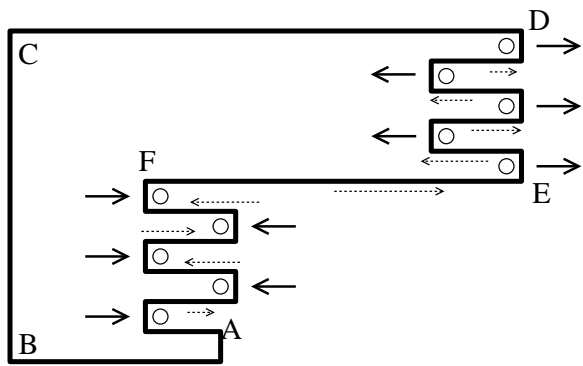}
  \caption{The pulley $C_h$, in case $h=5$}
  \label{fig: pulley}
\end{figure}

Consider the figure \vref{fig: pulley}.  The thick line $ABCDEF$ is
the curve $\theta \mapsto C_5(\theta,0)$. The thick arrows represent
the normal part $\pi_N \derpar v C_5$ of the velocity, while the thin
dashed arrows represent the tangent part $\pi_T \derpar v C_5$ of the
velocity.  The circles are just for fun, and represent the wheels of
the pulley.

The movement of $C_5(\theta,v)$, that is, its evolution in $v$, is so
described: the part $ABCD$ is still, that is, it is constant in $v$;
in the part $DE$ (respectively $FA$) of the curve, vertical segments move 
apart (resp., together)  as the thick arrows indicate, 
with horizontal velocity with norm $|\pi_N \derpar v C_5|=1$;
as the curve unravels, it is forced to move also parallel to itself.

The generic curve $C_h$ has $2h$ wheels: $h$ wheels in section $DE$,
to pull apart, and $h$ wheels in section $AF$, to pull together;
the horizontal tracts in $AF$ and $DF$ are of lenght
$\propto 1/h$, so the tract $AF$ straightens up in a time
$\Delta v\propto 1/h$: at that moment, the movement inverts:
so while $v\in[0,1]$, the cycle repeats for $h$ times.

In this case, the tangent velocity in section $DF$ is $h$ times the
normal velocity in sections $AF$ and $DE$. Then, if we choose the
normal velocity to have norm 1, the tangent velocity will explode
when $h\to\infinity$.
This means that the family of homotopies $C_h$
will not be compact in $H^1(I\to\real^1)$.

\bigskip

The first objection that comes to mind when reading
the above is ``\emph{this example is not showing any problem with
  the curve itself, it is just giving problems with the
  parameterization of the curve}''.

Indeed we may reparameterize the curves so that
the tangent velocity will \underline{not} explode
when $h\to\infinity$:  by using  \ref{prop: black magic},
we obtain that $\pi_{\tilde T} \derpar v \tilde C=0$.

Remembering remark \ref{rem: not both}, we understand that
this is not going to help, though. So we point out this other problem.
\begin{center}\em
  Let $\lambda=\lambda(v,h)$ be the distance from feature point $D$ to feature
  point $E$. Then $\derpar v \lambda\to \infinity$ if $h\to\infinity$. 
\end{center}

\section{Conformal metrics}
\label{sec: conformal}

\begin{figure}[ht]
  \centerline{
    \hfill
    \begin{minipage}[c]{.4\textwidth}
      \psfig{figure=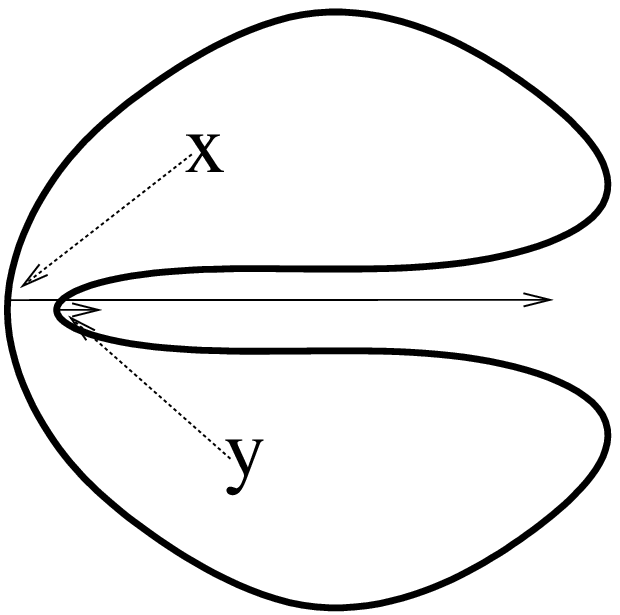,width=.95\textwidth}
    \end{minipage}
    \hfill
    \begin{minipage}[c]{.4\textwidth}
      \psfig{figure=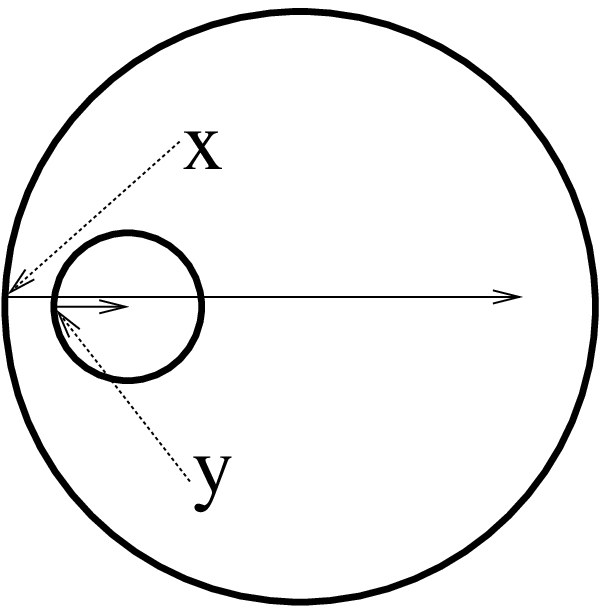,width=.95\textwidth}
    \end{minipage}
    \hfill
  }
  \caption[Intersections induced by flow.]{\label{fig:embedding}
    Intersections induced by flow \eqref{eq:heatflow-MM}: for big
    choices of $A$, the point $x$ will travel faster than the point $y$
    (and in the same direction) and will eventually cross it.
  }
\end{figure}

We recall at this point that the well known {\em geometric heat flow}
($C_t=C_{ss}$)%
\footnote{Recall that we write $C_v$ for $\derpar v C$ (and so on),
  to simplify the derivations
}
is truly the {\em gradient descent} for the Euclidean
arclength of a curve with respect to the $H^0$ metric.
Unfortunately, given the pathologies encountered thus far with $H^0$, we
see that this famous flow is not a gradient flow with respect to
a well behaved Riemannian metric. If we propose a different metric,
the new gradient descent flow for the Euclidean arclength of a curve
will of course be entirely different.
For example, the metric \eqref{eq:G MM} proposed
by \cite{Michor-Mumford} yields the following gradient flow for arclength:
\begin{equation} \label{eq:heatflow-MM}
  C_t = \frac{C_{ss}}{1+A \, C_{ss}\cdot C_{ss}}
      = \frac{\kappa}{1+A\kappa^2}N
\end{equation}
Notice that the normal speed in \eqref{eq:heatflow-MM} is {\em not
monotic in the curvature}; and therefore the flow
\eqref{eq:heatflow-MM} will not
share the nice properties of the  {\em geometric heat flow}
($\partial_t C=\partial_{ss} C$). For example, embedded curves
do not always remain embedded under this new flow, as illustrated
in Fig.~\vref{fig:embedding}.

Given the pathologies of $H^0$ we have no choice but to propose a new
metric if we wish to construct a well behaved Riemannian geometry on
the space of curves. However, we may seek a new metric whose gradient
structure is as similar as possible to that of the $H^0$ metric. In
particular, for any functional $E:M\to\real$ we may ask that the
gradient flow of $E$ with respect to our new metric be related to
that gradient flow of $E$ with respect to $H^0$ by only a time
reparameterization. In other words, if $C(t)$ represents a gradient
flow trajectory according to $H^0$ and if $\hat{C}(t)$ represents the
gradient flow trajectory according to our proposed new metric, then
we wish that
$$\hat{C}(t)=C(f(t))$$
for some positive time reparameterization $f:\real\to\real$, $\dot{f}>0$.
The resulting gradient flows will then be related as follows.
\begin{equation}\label{eq:reparam}
  \hat{C}_t= \dot{f}(t)\,C_t
\end{equation}

The only class of new metrics that will satisfy (\ref{eq:reparam}) are
{\em conformal} modifications of the original $H^0$ metric, which we
will denote by $H^0_\phi$. Such metrics
are completely defined by combining the original $H^0$ metric with a
positive {\em conformal factor} $\phi:M\to\real$ where $\phi(c)>0$ may
depend upon the curve $c$. The relationship between the inner products
is given as follows.
\begin{equation}\label{eq:conf_def}
  \Big\langle h_1,h_2\Big\rangle_{\!H^0_\phi}=
  \phi(c)\,\Big\langle h_1,h_2\Big\rangle_{\!H^0}
\end{equation}
Note that for any energy functional $E$ of curves $C(t)$ we have
the following equivalent expressions, where the first and last expressions
are by definition of the gradient and the middle expression comes from
the definition (\ref{eq:conf_def}) of a conformal metric.
\begin{equation}\label{eq:conf_deriv}
  \deriv{t}E(C(t))
   =\left\langle \pderiv[C]{t},
      \underbrace{\nabla^\phi E(C)}_
        {\parbox{3.5em}{\vspace{-2.5ex}\tiny\begin{center}Conformal\\Gradient\end{center}}}
    \right\rangle_{\!\!\raisebox{4ex}{$\scriptstyle H^0_\phi$}}
   \!\!\!\!\!
   =\phi\left\langle \pderiv[C]{t},
      \underbrace{\nabla^\phi E(C)}_
        {\parbox{3.5em}{\vspace{-2.5ex}\tiny\begin{center}Conformal\\Gradient\end{center}}}
    \right\rangle_{\!\!\raisebox{4ex}{$\scriptstyle H^0$}}
   \!\!\!\!\!
   =\left\langle \pderiv[C]{t},
      \underbrace{\nabla E(C)}_
        {\parbox{3.5em}{\vspace{-2.5ex}\tiny\begin{center}Original\\Gradient\end{center}}}
    \right\rangle_{\!\!\raisebox{4ex}{$\scriptstyle H^0$}}
\end{equation}
We see from (\ref{eq:conf_deriv}) that the 
{\em conformal gradient differs only in magnitude} from the original
$H^0$ gradient
$$
  \nabla^\phi E = \frac 1{\phi} \nabla E
$$
and therefore the {\em conformal gradient flow differs only in speed}
from the $H^0$ gradient flow.
$$ \pderiv[C]{t} = - \nabla^\phi E(C) = -\frac 1{\phi(C)} \nabla E(C) $$
As such and as we desired, the
{\em solution differs only by a time reparameterization $f$} given by
$$ \dot{f}=\frac 1{\phi(C)}$$

The obvious question now is how to choose the conformal factor.

A first suggestion is in the following theorem
\begin{Proposition}\label{prop: conforme non deg}
\footnote{We thanks prof. Mumford for suggesting this result.}
  Suppose that 
  \begin{equation}
    \min_c(\phi(c)/\len(c))=a>0\label{eq:phinondeg}
  \end{equation}
  Consider an homotopy $C$ connecting two curves 
  $c_0=C(\cdot,0),c_1=C(\cdot,1)$,
  and its $H^0_\phi$--energy
  \[\int_0^1 \phi(C(\cdot,v))
  \int_{S^1}|\pi_N \derpar v C|^2 |\dot C| \, d\theta \,dv \] 

  Up to  reparameterization \ref{prop: by arc par},
  $|\dot C(\theta,v)|=\len(C(\cdot,v))/2\pi$ so we can
  rewrite the energy (using \eqref{eq:|VxW|}) as
  \begin{eqnarray*}
    \int_0^1\frac{ 2\pi \phi(C(\cdot,v))}{\len(C(\cdot,v))}
    \int_{S^1}|\derpar v C\times\dot C |^2 \, d\theta \,dv \ge 2\pi a
    \int_0^1 \int_{S^1}|\derpar v C\times\dot C |^2 \, d\theta \,dv \ge
    \\  \ge a \left(\int_0^1 \int_{S^1}|\derpar v C\times\dot C | \,
      d\theta \,dv \right)^2
  \end{eqnarray*}    
  the rightmost term is the square of the area swept by the
  homotopy.

  So, if $c_0\neq c_1$, and there do not exist a homotopy
  connecting $c_0,c_1$ with zero area, then $d(c_1,c_2)>0$. 
\end{Proposition}

Although
we already know that the $H^0$ metric is not very useful, we may obtain
a lot of insight into how to choose the conformal factor $\phi$ by
observing the structure of the minimizing flow (which turns out to
be unstable) for the $H^0$ energy in the space of homotopies. We will
then try to choose the conformal factor in order to counteract the
unstable elements of the $H^0$ flow.

\subsection{The Unstable $H^0$ Flow}

\subsubsection{Geometric parameters $s$ and $\v$}
\label{sec: geo par}

We have denoted by $u\in[0,1]$ a parameter which traces out
each curve in a parameterized homotopy $C(u,v)$ and we have denoted
by $v\in[0,1]$ the parameter which moves us from curve to curve
along the homotopy. Note that both of these parameters are arbitrary
and not unique to the geometry of the curves comprising the homotopy.
We now wish to construct more geometric parameters for the homotopy
which will yield a more meaningful and intuitive expression for the
minimizing flow we are about to derive. The most natural substitute
for the curve parameter $u$ is the arclength parmeter $s$. We must
also address the parameter $v$, however. While
$v$ as a parameter ranging from 0 to 1 seems to have little to do
with the arbritrary choice of the curve parameter $u$,
the differential operator $\pderiv{v}$ depends heavily upon this
prior choice. The desired effect of differentiating along the homotopy
is mixed with the undesired effect of differentiating along the contour
if flowing along corresponding values of $u$ between curves in the
homotopy requires some motion along the tangent direction.
To see the dependence of $\pderiv{v}$ on $u$, note that $C(u,v)$ and
$\hat{C}(u,v)$ where
$$\hat{C}(u,v)=C\left(u^{(1+v)},v\right)$$
constitute the same homotopy geometrically,
and yet $\pderiv[C]{v}\ne\pderiv[\hat{C}]{v}$.

We will therefore introduce the more geometric parameter $\v$ whose
corresponding differential operator $\pderiv{\v}$ yields the most
efficient transport from one curve to another curve along the homotopy
regardless of ``correspondence'' between values of the curve parameters.
It is clear that such a transport must always move in the normal direction
to the underlying curve since tangential motion along any curve does not
contribute to movement along the homotopy. More preceisely, we define
the parameteres $s$ and $\v$ in terms of $u$ and $v$ as follows.
$$
\pderiv{s}=\frac{1}{\|C_u\|}\pderiv{u} \qquad\mbox{and}\qquad
\pderiv{\v}=\pderiv{v}-\big(C_v\cdot C_s\big) \pderiv{s}
$$

\subsubsection{$H^0$ Minimizing Flow}

Suppose we now consider a time varying family of
homotopies $C(u,v,t):[0,1]\times [0,1]\times (0,\infty)\to \real ^n$
and compute the derivative of the $H^0$ energy along
this family. Note that the $H^0$ energy, in terms of
the new parameters $s$ and $\v$ may be simply expressed
as follows (since $\pi_N C_{\v}=C_{\v}$).
\begin{equation}
  E(t)=\int_0^1\int_0^L \big\|C_{\v}\big\|^2\,ds\,dv
\end{equation}
In the appendix, we show that the derivative of $E$ may be
expressed as follows.
\begin{equation}
  E'(t) =
      -2\int_0^1\int_0^L
        C_t\cdot\Big(
          C_{\v\v}\!-(C_{\v\v}\!\!\cdot C_s)C_s
         -(C_{\v}\!\!\cdot C_{ss})C_{\v}
         +\frac{1}{2}\|C_\v\|^2C_{ss}
        \Big)
      \,ds\,dv
\end{equation}
In the planar case, $C_{\v}$ and $C_{ss}$ are linearly dependent
(as both are orthogonal to $C_s$) which means that
\begin{equation}
  (C_{\v}\cdot C_{ss})C_{\v}=(C_{\v}\cdot C_{\v})C_{ss}=\|C_\v\|^2C_{ss}
\end{equation}
and therefore
\begin{equation}
E'(t)=
  -2\int_0^1\int_0^L
    C_t\cdot\Big(
      \big(C_{\v\v}-(C_{\v\v}\cdot C_s)C_s\big)-\frac{1}{2}\|C_\v\|^2C_{ss}
    \Big)
  \,ds\,dv
\end{equation}
by which we derive the minimization flow 
\[  C_t = C_{\v\v}-(C_{\v\v}\cdot C_s)C_s - \frac{1}{2}\|C_\v\|^2C_{ss} \]
which is geometrically equivalent to the following more simpler
flow (by adding a tangential component):
\begin{equation} \label{eq:flowH0}
  C_t = C_{\v\v} - \frac{1}{2}\|C_\v\|^2C_{ss}
\end{equation}

Note that the flow (\ref{eq:flowH0}) consists of two orthogonal
diffusion terms. The first term $C_{\v\v}$ is stable as it
represents a forward diffusion along the homotopy,
while the second term $-\|C_\v\|^2C_{ss}$ is an unstable
backward diffusion term along each curve.
Indeed, numerical experiments show a behaviour that parallels 
the phenomenon described in \S\ref{sec: tesselation}.

\subsection{Conformal Versions of $H^0$}

We now define the conformal $H^0_\phi$ energy (when the conformal
factor $\phi$ is a function of the arclength $L$ of each curve) as
\begin{equation}
  E_{\phi}(t)=\int_0^1 \phi(L) \int_0^L \big\|C_{\v}\big\|^2\,ds\,dv
\end{equation}
Once again we compute (in the appendix) the derivative of this
energy along a time varying family of homotopies $C(u,v,t)$.

\begin{eqnarray*}
  E_\phi'(t) =
      -\int_0^1\int_0^L
        C_t\cdot
        \hspace{-2em}
        && \Big(
            2\phi'L_{\v} C_{\v}
           +2\phi C_{\v\v}-2\phi(C_{\v\v}\cdot C_s)C_s \\
        && -2\phi(C_{\v}\cdot C_{ss})C_{\v}
           +(\phi m+\phi'M)C_{ss}
          \Big)
      \,ds\,dv
\end{eqnarray*}
where
\[
  m=\|C_\v\|^2 \qquad\mbox{and}\qquad
  M=\int_0^L m \,ds = \int_0^L \|C_\v\|^2 \,ds.
\]
As before, we now consider the planar case in
which $C_{\v}$ and $C_{ss}$ are linearly dependent
and therefore $(C_{\v}\cdot C_{ss})C_{\v}=mC_{ss}$,
yielding
\begin{equation}
E'(t)= -2
   \int_0^1\int_0^L
    \!C_t\cdot\Big(
      \phi\big(C_{\v\v}\!\!-(C_{\v\v}\!\!\cdot C_s)C_s\big)
     +\phi'L_{\v} C_{\v}
     +\frac{1}{2}(\phi'M-\phi m)C_{ss}
    \Big)
  \,ds\,dv
\end{equation}
from which we obtain the following minimizing flow
\[ C_t =  
      \phi\big(C_{\v\v}-(C_{\v\v}\cdot C_s)C_s\big)
     +\phi'L_{\v} C_{\v}
     +\frac{1}{2}(\phi'M-\phi m)C_{ss}
\]
which is geometrically equivalent (by adding a tangential
term) to
\begin{equation}\label{eq:flowConf}
   C_t =  
      \phi C_{\v\v}
     +\phi'L_{\v} C_{\v}
     +\frac{1}{2}(\phi'M-\phi m)C_{ss}
\end{equation}

\subsubsection{Stable conformal factor}

To stabilize the flow in last equation, we look for a $\phi$ such that

\begin{equation}\label{eq: phi stab}
  \phi'M-\phi m \ge 0 \qquad\mbox{for all }(s,\v)
\end{equation}
or (assuming $M\ne 0$)
\begin{equation}
  \frac{\phi'}{\phi}=(\log\phi)'
  \ge \frac{m}{M} \qquad\mbox{for all }(s,\v)
\end{equation}
One way to satisfy this is to choose
\begin{equation}\label{eq: lambda stab}
 (\log\phi)'=\max_{s,\v}\frac{m}{M}\doteq\lambda
\end{equation}
giving us
\begin{equation}
 \phi=e^{\lambda L}
\end{equation}
yielding the following flow of homotopies
\begin{equation}
 C_t= e^{\lambda L}\Big( 
        2C_{\v\v}+2\lambda L_{\v}C_{\v}+(\lambda M-m)C_{ss}
      \Big)
\end{equation}

The choice of having $\phi=e^{\lambda L}$ 
satisfies \eqref{eq:phinondeg}; and it agrees
also with the discussion in \S\ref{sec: exa conf lsc}, 
that hints that the energy $E(C)$ associated to the $H^0_\phi$ metric
may be lower semi continuous when 
$\phi(c)\ge \len(c)$.
The above conformal metric does not entail an unique Riemannian Metric 
on the space of curves: indeed the choice of $\lambda$ depends on
the homotopy itself.

In the numerical experiments shown in this paper, we chose $\lambda$ 
to satisfy \eqref{eq: lambda stab} at time $t=0$ and found that
this was enough to stabilize the flow up to convergence. However,
we have no mathematical proof of this phenomenon.

\subsection{Numerical results}
\label{sec: numerical}
We note that the minimizing flow (\ref{eq:flowConf}) consists of
two stable diffusion terms and a transport term. As such, we have
the option to utilize level set methods in the implementation of
(\ref{eq:flowConf}).

We represent the evolving homotopy $C(u,v,t)$ as an evolving surface
$S(u,v,t)$
  $$S(u,v,t)=(C(u,v),v,t)$$
We then perform a {\em Level Set Embedding} of this surface into a 4D scalar
function $\psi$ such that
  $$\psi\big(C(u,v,t),v,t\big)=0.$$

The goal is now to determine an evolution for $\psi$ which yields
the evoluton (\ref{eq:flowConf}) for the level sets of each of
its 2D cross-sections. Differentiating
$$
  \frac{d}{dt}\Big(\psi\big(x(u,v,t),y(u,v,t),v,t\big)=0\Big)
  \quad\longrightarrow\quad
  \psi_t+\nabla\psi\cdot C_t=0
$$
where $\nabla\psi=(\psi_x,\psi_y)$ denotes the 2D spatial gradient of each
2D cross-section of $\psi$, and substituting (\ref{eq:flowConf}),
noting that $N=\nabla\psi/\|\nabla\psi\|$, yields the corresponding
{\em Level Set Evolution}.
  \begin{eqnarray*}
  \psi_t\!=&& \hspace{-7mm}
     \psi_{vv}
    \!-\!\frac{2\psi_v}{\|\nabla\psi\|^2}(\nabla\psi_v\!\cdot\!\nabla\psi)
    \!+\!\frac{\psi_v^2}{\|\nabla\psi\|^4}
    \big(\nabla^2\psi\nabla\psi\big)\!\cdot\!\nabla\psi
  \\ && \hspace{-7mm} \nonumber
    -\frac{1}{2}\left(
      \frac{\psi_v^2}{\|\nabla\psi\|^2}
     -\lambda\int_0^L\frac{\psi_v^2}{\|\nabla\psi\|^2}\,ds
    \right)
    \nabla\cdot\left(\frac{\nabla\psi}{\|\nabla\psi\|}\right)\|\nabla\psi\|
   \!+\!\lambda L_v\psi_v
  \end{eqnarray*}
Note that for simplicity we have dropped the factor
$e^{\lambda L}$ from (\ref{eq:flowConf}) since we are guaranteed
that this factor is always positive. As a result, we do not change
the steady-state of the flow by omitting this factor.

\begin{figure}[htbp]
  \centerline{
    \begin{tabular}{cc}
      \psfig{figure=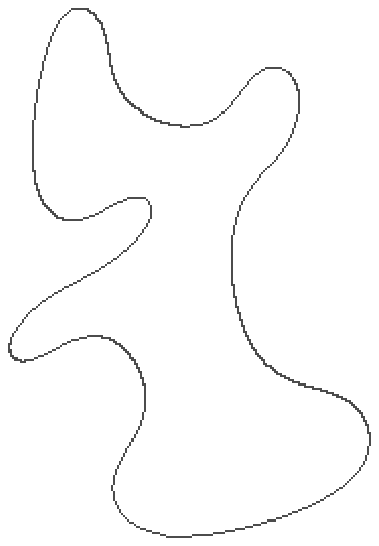,width=2in} &
      \psfig{figure=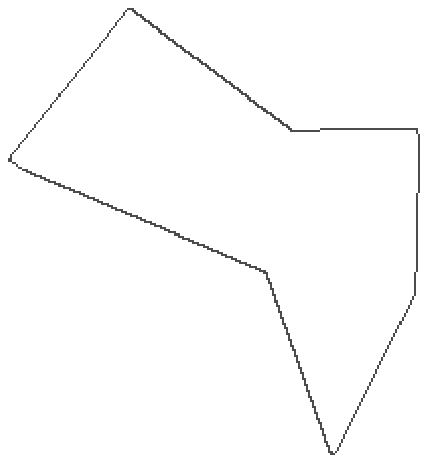,width=2in} \\
      Curve $c_0$ & Curve $c_1$ \\
    \end{tabular}
  }
  \caption{Curves  $c_0$ and $c_1$}
  \label{fig: num c0 c1}
\end{figure}

\begin{figure}[htbp]
   \psfig{figure=Figures/contour00,width=1.09in}\hspace{-.25in}
   \psfig{figure=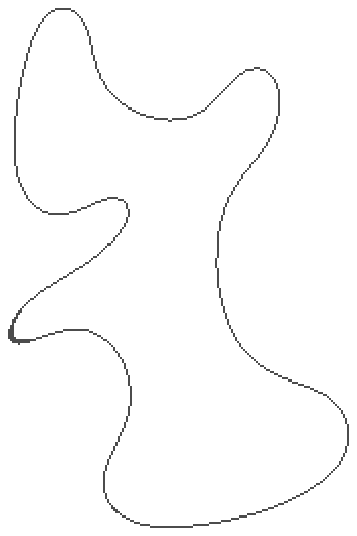,width=1.09in}\hspace{-.25in}
   \psfig{figure=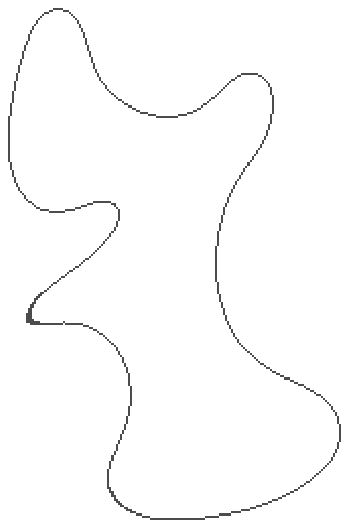,width=1.09in}\hspace{-.25in}
   \psfig{figure=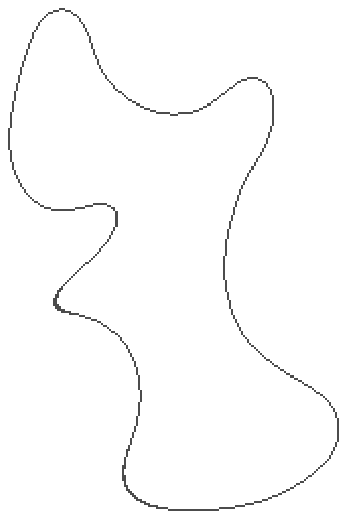,width=1.09in}\hspace{-.25in}
   \psfig{figure=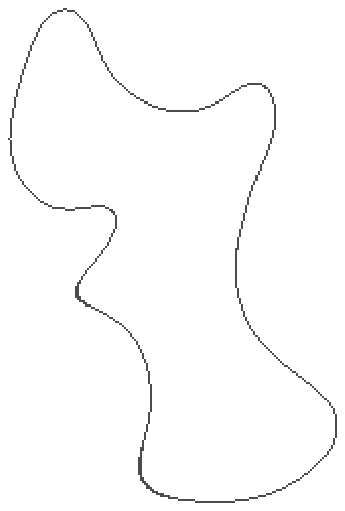,width=1.09in}

   \psfig{figure=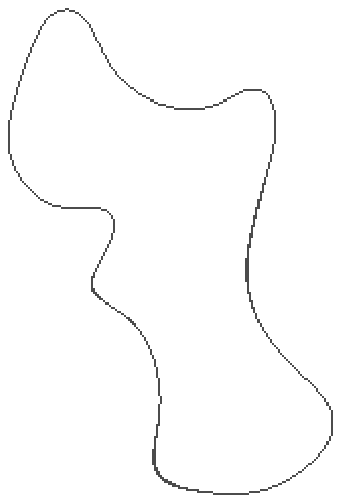,width=1.09in}\hspace{-.25in}
   \psfig{figure=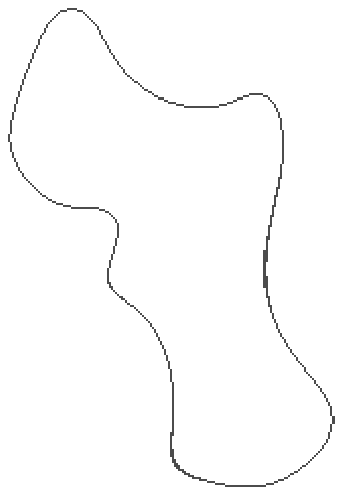,width=1.09in}\hspace{-.25in}
   \psfig{figure=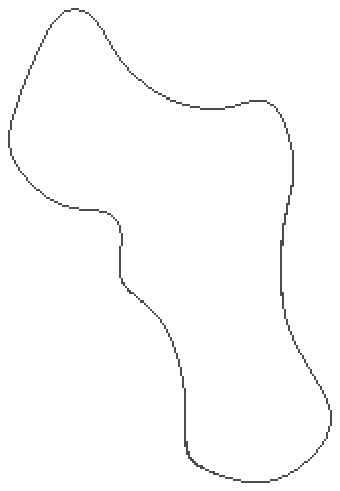,width=1.09in}\hspace{-.25in}
   \psfig{figure=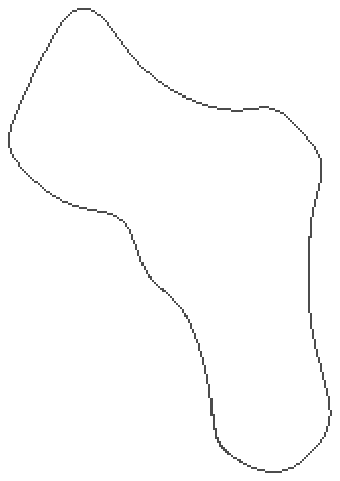,width=1.09in}\hspace{-.25in}
   \psfig{figure=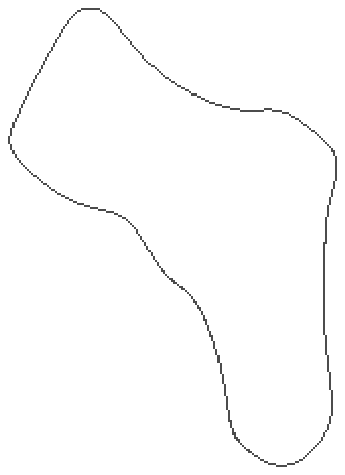,width=1.09in}

   \psfig{figure=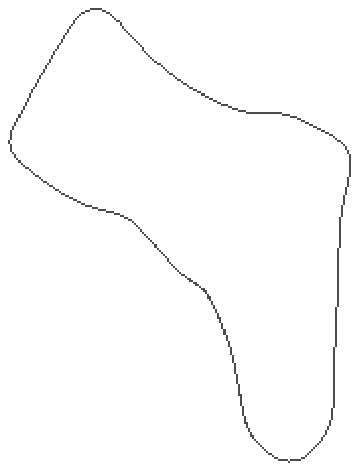,width=1.09in}\hspace{-.25in}
   \psfig{figure=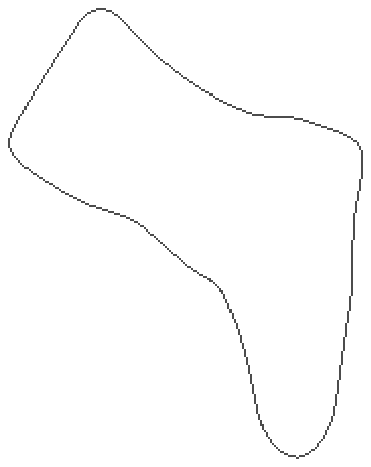,width=1.09in}\hspace{-.25in}
   \psfig{figure=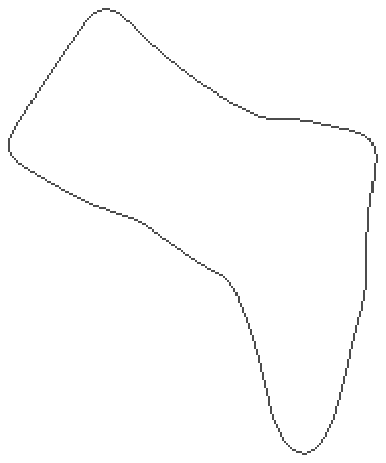,width=1.09in}\hspace{-.25in}
   \psfig{figure=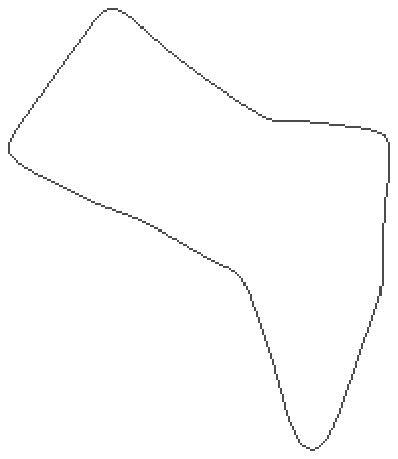,width=1.09in}\hspace{-.25in}
   \psfig{figure=Figures/contour63,width=1.09in}\centering
   
   \caption{Slices of homotopy}
   \label{fig: num slice}
 \end{figure}

\begin{figure}[htbp]
  \psfig{figure=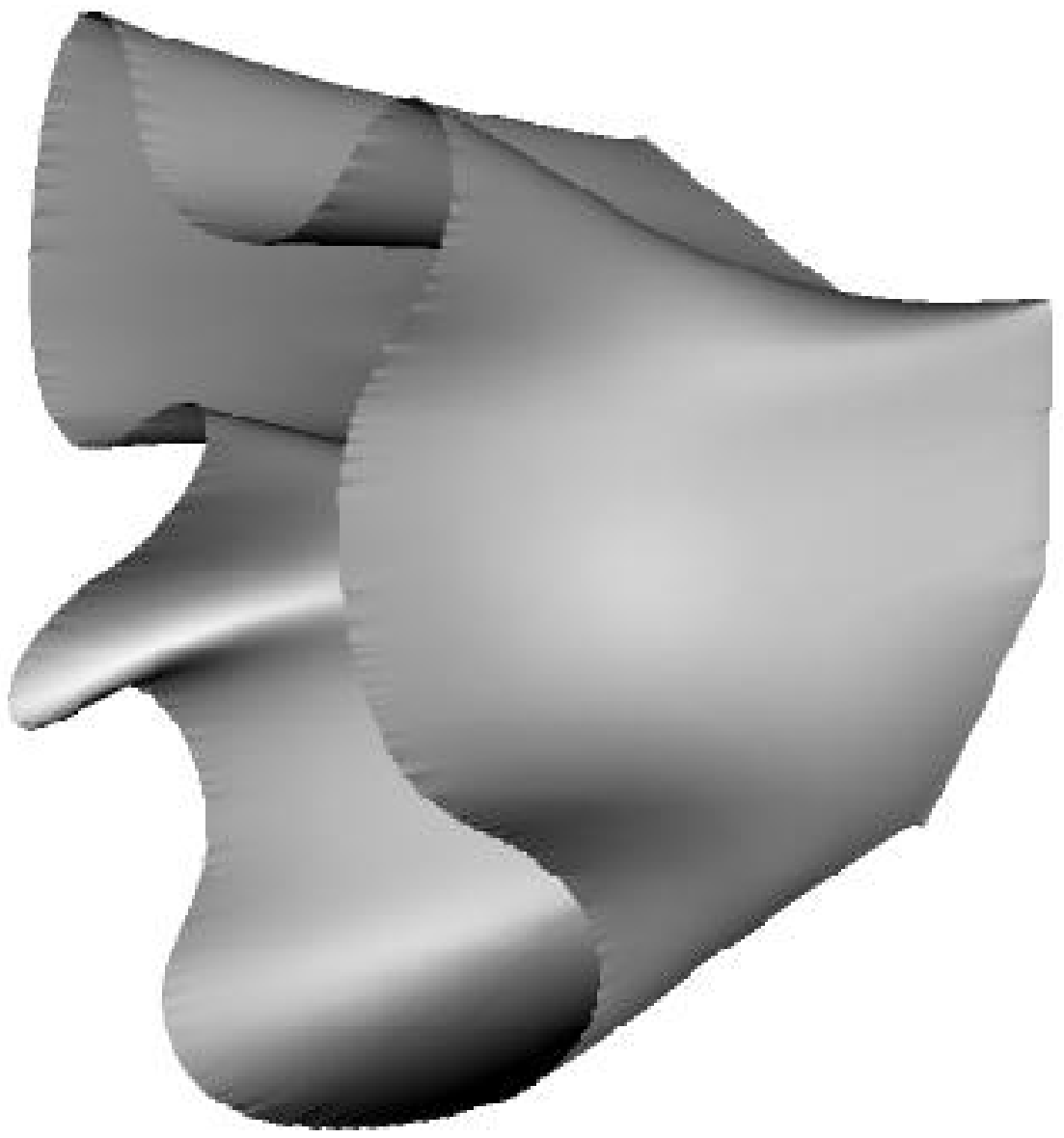,width=.495\textwidth}
  \psfig{figure=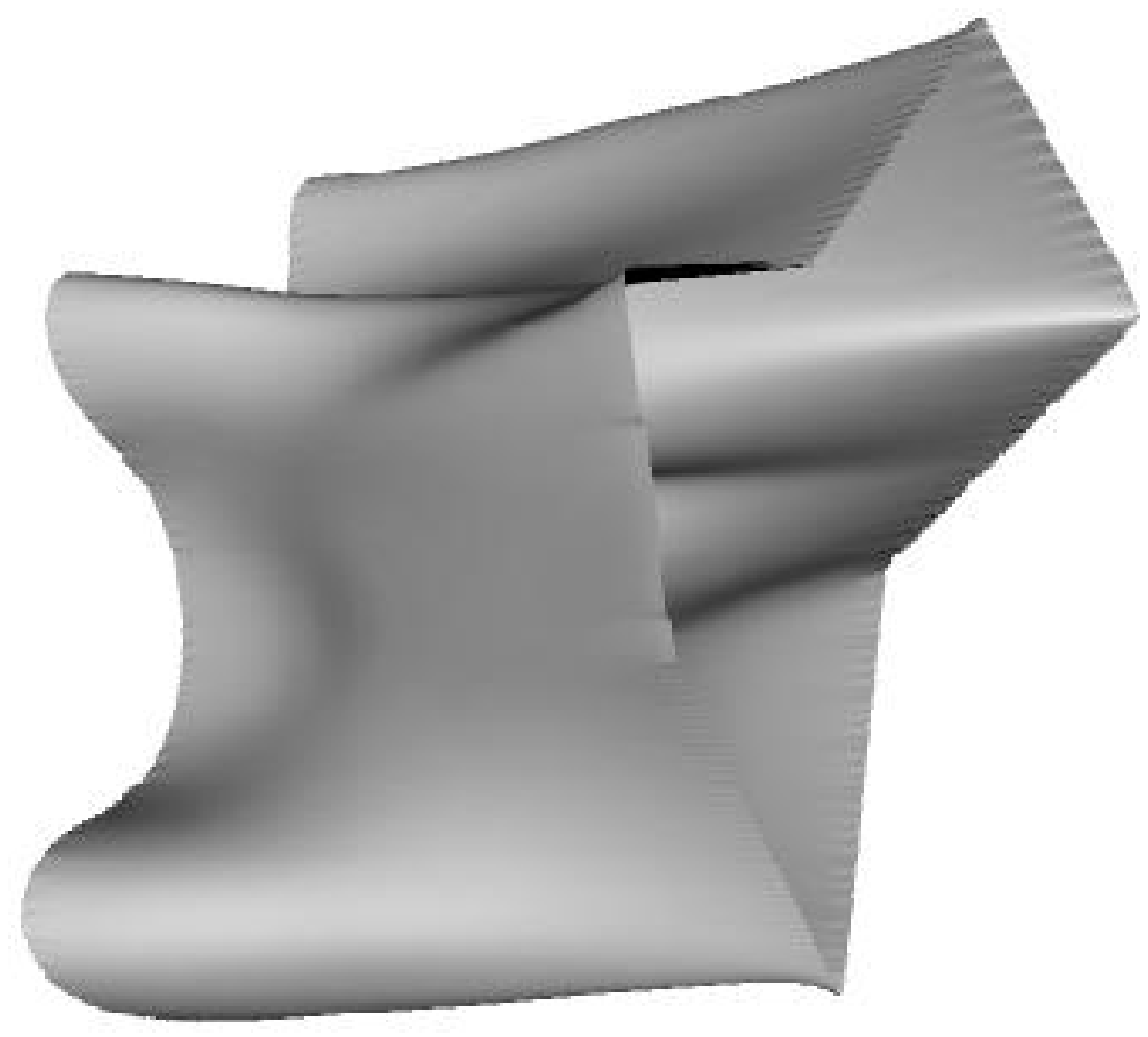,width=.495\textwidth}
   \caption{Surface of homotopy}
   \label{fig: num surf}
\end{figure}

If we numerically compute the geodesic of two curves $c_0,c_1$ in figure
\vref{fig: num c0 c1},
we obtain the geodesic, which is represented, by slicing it in
figure \vref{fig: num slice} and as a surface in figure \vref{fig: num surf}.

\subsection{Example}
\label{sec: exa conf lsc}
\begingroup
\def\t{u}
\def\e{\varepsilon}
\def\l{\lambda}
\def\a{\alpha}
This example shows why we think that the 
conformal energy $E(C)$ may be lower semicontinuous (on planar curves)
in the case when  $\phi(c)\ge \len (c)$.

Fix $0< \e<1/2$ and $\l\ge 0$.

Suppose 
\begin{equation}
  c(\t) = \oldcases {
 (\t,\t\l)  & if $\t\in [0,\e]$ \cr
 (\t,(2\e-\t)\l)  & if $\t\in [\e,2\e]$ \cr
 (\t,0)  & if $\t\in [2\e,1]$ }
\label{eq:kjedve}
\end{equation} 
is the curve in figure \vref{fig: c}.

\begin{figure}[htbp]
  \centering
  \includegraphics{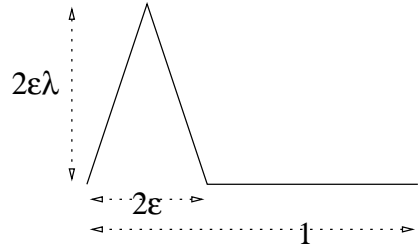}
  \caption{The curve from eq.~\eqref{eq:kjedve}}
  \label{fig: c}
\end{figure}

We define the homotopy $C:[0,1]\times[0,1]\to \real^2$ by
\[\tilde C(\t,v) = c(\t)+ (0,v)\]

Let with $C(\t,v)=(\t,v)$ be the identity.

We may tesselate, as explained in point
\ref{item:tessel} in \ref{exa: E not lsc}, to
build a sequence of homotopies $C_h$ such that 
$C_h\to_h C $ in $L^\infinity$ and 
\[\derpar{\t,v}C_h\wto^*_h \derpar{\t,v}C \mbox{ weakly* in
$L^\infinity$},
\]
 and $E(C_h)=E(C)$.

Now we compute
\[ \len(c) = 1+2\e(\sqrt{1+\l^2}-1) = 1+ 2\e \a\]
where we define $\a=\sqrt{1+\l^2}-1$ for convenience.
Note that $\a \ge 0$.

We  compute the energy (using the identity \eqref{eq:|VxW|})
\begin{eqnarray*}
  \label{eq:E ex conforme}
  E(C_h)=E(\tilde C) =
  \int_0^1\int_0^1 \frac{ \phi(C)}{|\derpar \t C|}d\t dv
  =  \\  =
  ( 1-2\e)  \phi(C)  +2\e  \frac{ \phi(C)}{ \sqrt{1+\l^2} } 
  =\phi(C) \left(  1-2\e+2\e\frac 1{\a+1}  \right )
  = \\ = 
  \phi(C)\left(  1-2\e\frac{\a} {\a+1}   \right )
  \ge \\ \ge
  (1+2\e\a)\left(  1-2\e\frac{\a} {\a+1}   \right )=
  1 +2\e\left(\a - \frac{\a} {\a+1} -\frac{2\e \a^2} {\a+1}  \right )
  = \\  =  1 +2\e \frac{ \a^2-2\e\a^2}{\a+1}= 
  1 +2\e\a^2 \frac{ 1-2\e}{\a+1} \ge 1 = E^N(C)
\end{eqnarray*}
\endgroup

\appendix

\section{More on Finsler metrics} \label{app:finsler}

We now provide two more results
on Finsler metrics \ref{def: Finsler}, for convenience of the reader.
The first result explains the relationship between the
length functional $\len(\xi)$ and the energy functional $E(\xi)$
\begin{Proposition}\label{prop: len and E}
    \def\A{{\mathcal A}} Fix $x,y$ in the following, and let $\A$ be
    the class of all locally Lipschitz paths $\gamma:[0,1]\to M$
    connecting $x$ to $y$.
    
  These are known properties of the the length and the energy.
  \begin{itemize}
  \item If $\xi,\gamma$ are locally Lipschitz and $\phi$ is a monotone
    continuous function such that $\xi=\gamma \circ \phi$
    and $\phi(0)=0,\phi(1)=1$ then
    $\Len\gamma=\Len\xi$.
    
  \item 
    In general, by H\"older inequality, $E(\gamma)\ge \Len(\gamma)^2$.
    
  \item 
    If $\gamma$ provides a minimum of $\min_\A {E(\gamma)}$, then it
    is also a minimum of $\min_\A {\Len(\gamma)}$ in the same class,
    ${E(\gamma)} = \Len(\gamma)^2$, and moreover $|\dot
    \gamma(t)|_{\gamma(t)}$ is constant in $t$, that is, $\gamma$ has
    constant velocity.
    
  \item 
    If $\xi$ provides a minimum of $\min_\A {\Len(\gamma)}$, then
    there exists a monotone continuous function $\phi$ and a path
    $\gamma$ such that $\xi=\gamma \circ \phi$, and $\gamma$ is a
    minimum of $\min_\A {E(\gamma)}$.  ~\footnote{in writing
      $\xi=\gamma \circ \phi$, we, in a sense, define $\gamma$ by a
      \emph{pullback} of $\xi$: see 2.29 in \cite{ACM:Asym}.  Note
      that we could not write, in general, $\gamma=\xi \circ
      \phi^{-1}$: indeed, it is possible that a minimum of $\min_\A
      {\Len(\gamma)}$ may stay still for an interval of time; that is,
      we must allow for the case when $\phi$ is not invertible.}

\end{itemize}
  \begin{proof}
    By 2.27, 2.42, 3.8, 3.9 in \cite{ACM:Asym} and
    4.2.1 in \cite{Ambrosio-Tilli}.
  \end{proof}
\end{Proposition}

\medskip

This second result is the version of the Hopf--Rinow theorem
\ref{Hopf-Rinow} to the case of generic metric spaces.
\begin{Theorem}[Hopf-Rinow]\label{metric Hopf-Rinow}
   Suppose that the metric space
   $(M,d)$ is locally compact and path-metric, then these are equivalent:
   \begin{itemize}
   \item the metric space $(M,d)$ is complete,
   \item closed bounded sets are compact
 \end{itemize}
 and both imply that any two points can be connected by a minimal length
 geodesic.
\end{Theorem}
A proof is in \S1.11 and \S1.12 in Gromov's  \cite{Gromov:met},
or in Theorem~1.2 in   \cite{ACM:Asym} (which holds
also in the asymmetric case).

\section{Proofs of \S{\protect\ref{SEC:SRIV}}}
\label{sec: proofs Sriva}

Let $Z$ be the set of all $\theta\in L^2([0,2\pi])$ such that
$\theta(s)=a+k(s) \pi$ where $k(s)\in \integer$ is
measurable, $(a=2\pi-\int k)$,
and \[|\{k(s)=0 \mod 2\} | = |\{k(s)=1 \mod 2\} | = \pi\]
$Z$ is closed (by thm. 4.9 in \cite{Brezis}).
We see that $Z$ contains the (representations $\theta$ of) flat curves
$\xi$, that is, curves $\xi$ whose image is contained in a line;
one such curve is 
\[  \xi_1(s)= \xi_2(s) = 
\begin{cases}
  s/\sqrt 2 & s\in[0,\pi] \\
  (2\pi-s)/\sqrt 2 & s\in(\pi,2\pi]
\end{cases}
,\qquad
\theta= 
\begin{cases}
  \pi/2 & s\in[0,\pi] \\
  3\pi/2& s\in(\pi,2\pi]
\end{cases}
\]

We provide here the proof of
\ref{prop: Sriva M}: $M\setminus Z$ is a manifold.
\begin{proof}
  Indeed, suppose by contradiction that
  $\nabla \phi_1  , \nabla \phi_2, \nabla \phi_3$
  are linear dependant at $\theta\in M$, 
  that is, there exists $a\in\real^3,a\neq 0$ s.t.
  \[a_1 \cos(\theta(s))+a_2 \sin(\theta(s)) + a_3 = 0\]
  for almost all $s$; then, by integrating, $a_3 = 0$,
  therefore $a_1 \cos(\theta(s))+a_2 \sin(\theta(s)) = 0$ that means
  that $\theta\in Z$.
  See also \S3.1 in  \cite{Srivastava:AnalPlan}.
\end{proof}

This is the proof of \ref{prop: no rotation Srivastava}:
\begin{proof}
  Fix $\theta_0\in M\setminus Z$. Let $T=T_{\theta_0} M$ be the
  tangent at $\theta_0$. $T$ is the vector space orthogonal to
  $\nabla \phi_i(\theta_0)$ for $i=1,2,3$.
  Let $e_i=e_i(s)\in L^2\cap C^\infinity_c$
  be near $\phi_i(\theta_0)$ in $L^2$, so that the map
  $(x,y):T\times\real^3\to L^2$
  \begin{equation}
    (x,y) \mapsto 
    \theta=\theta_0+x + \sum_{i=1}^3  e_i y_i 
    \label{eq:loc coord}      
  \end{equation}
  is an isomorphism.
  Let $M'$ be $M$ in these coordinates; by the Implicit Function
  Theorem (5.9 in \cite{Lang:FDG}), there exists an open set
  $U'\subset T$, $0\in U'$, an open $V'\subset \real^3$, $0\in
  V'$, and a smooth function $f:U\to\real^3$ such that the local
  part $M'\cap(U'\times V')$ of the manifold $M'$ is the graph of
  $y=f(x)$.
    
  We immediatly define a smooth projection $\pi:U'\times V'\to M'$
  by setting $\pi'(x,y) = (x, f(x))$; this may be expressed in
  $L^2$; let $(x(\theta),y(\theta))$ be the inverse of
  \eqref{eq:loc coord} and $U=x^{-1}(U')$; we define the
  projection $\pi:U\to M$ by setting
  \[\pi(\theta) = \theta_0+x + \sum_{i=1}^3  e_i f_i(x(\theta)) \] Then
  \begin{equation}
    \pi(\theta)(s)-\theta(s) = \sum_{i=1}^3 e_i(s) a_i ~~,~~
    a_i\defeq (f_i(x(\theta))-y_i)\in\real \label{eq:pi(x)-x}    
  \end{equation}
  so if $\theta(s)$ is smooth, then $\pi(\theta)(s)$ is smooth.
  
  Let $\theta_n$ be smooth functions such that $\theta_n\to
  \theta$ in $L^2$, then $\pi(\theta_n)\to \theta_0$; if we choose
  them to satisfy $\theta_n(2\pi)-\theta_n(0)=2\pi h$, then, by
  the formula \eqref{eq:pi(x)-x}, $
  \pi(\theta)(2\pi)-\pi(\theta)(0) =2\pi h$ so that
  $\pi(\theta_n)\in M$ and it represents a smooth curve with the
  assigned rotation index $h$.
\end{proof}

\section{$E^N$ is ill-posed}
\label{sec: ill-posedness}

\begin{figure}[htbp]
  \begin{minipage}[tr]{.5\textwidth}
    \hfill 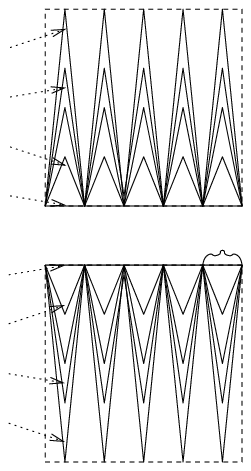
  \end{minipage}
  \begin{minipage}[tr]{.5\textwidth}
    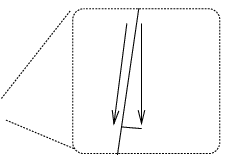
    \begin{eqnarray*}
      |\pi_N \derpar v C|^2 |\derpar \theta C|=
      \\ = |\derpar v C |^2|\derpar\theta C|\sin^2(\alpha)\sim
      \\ \sim \frac {|\derpar v C |^2}{ |\derpar\theta C|}
    \end{eqnarray*}
  \end{minipage}
  \caption{Artistic rendition of the homotopy $C_k$, from  \cite{Mumford:Gibbs}}
  \label{fig:zigzag}
\end{figure}

The result \ref{inf EN 0} following below
was inspired from a description 
of a similar phenomenon, found  on page 16 of the
slides  \cite{Mumford:Gibbs} of D. Mumford: 
it is possible to connect the two segments
$c_0(u)=(u,0)$ and $c_1(u)=(u,1)$ with a family of
homotopies $C_k:[0,1]\times[0,1]\to \real^2$ such that $E^N(C_k)\to_k 0$.
We represent the idea in figure \vref{fig:zigzag}.

We use the above idea to show that the distance induced
by $E^N$ is zero.
\footnote{
  We have recently discovered an identical proposition in
  \cite{Michor-Mumford}: we anyway propose this proof, since it is more
  detailed.
}

\begin{Proposition}[$E^N$ is ill-posed]
  \label{inf EN 0}
  \def\e{\epsilon}
  Fix $c_0$ and   $c_1$ to be two regular curves.
  We want to show that, $\forall \e>0$, there is a homotopy $C$
  connecting $c_0$ to $c_1$ such that $E^N(C)<\e$.

  \begin{enumerate}
  \item To start, suppose that $c_1$ is contained in the surface of a
    sphere, that is, $|c_1|=1$ is constant. Suppose also that
    $|\overdot c_1|=1$.
    Consider the linear
    interpolant, from the origin to $c_1$:
    $$C(\theta,v) = v c_1(\theta) $$
    The image of this homotopy is a cone.
    
    We want to play a bad trick to the linear interpolant: we define
    a homotopy whose image is the cone, but that moves points with
    different speeds and times.  Let $\e=\pi/k$ in the following;
    we define the sawtooth $Z:S^1 \to[0,\e]$
    $$Z(\theta) =\oldcases{ \theta & if $\theta\in[0,\e]$ \cr
      (2\e-\theta) & if $\theta\in[\e,2\e]$ \cr (\theta-2\e) & if
      $\theta\in[2\e,3\e]$ \cr (4\e-\theta) & if $\theta\in[3\e,4\e]$
      \cr \cdots }$$
    
    (note that $Z(\theta)+Z(\theta+\e)=\e$, 
    $Z(\theta)=Z(-\theta)$)

    Let
    $$C_k(\theta,v) = c_1(\theta) \frac{2v} {\e} Z(\theta)$$
    for
    $v\in[0,1/2]$, and
    \begin{eqnarray*}
      C_k(\theta,v) = c_1(\theta) \left(1-2\frac{1-v}{\e} Z(\theta+\e)\right)
    \end{eqnarray*}
    for $v\in[1/2,1]$.

    The energy $E(C_k)$ is splitted in two parts for $v$ and in
    $2k$ equal parts in $\theta$, so we compute the energy only for
    two regions, and then multiply by $2k$.

    \begin{itemize}
    \item In region $\theta\in[0,\e]$ $v\in[0,1/2]$, we have 
 
      $$C_k(\theta,v) = c_1(\theta) 2\frac 1 {\e} v \theta $$

      $$\derpar v C_k = c_1(\theta) 2\frac 1 {\e}\theta$$
      $$\derpar\theta C_k= \overdot c_1(\theta) 2\frac v {\e}\theta
      +c_1(\theta) 2\frac v {\e}$$
      and
      $$|\derpar\theta C_k|^2= 4 v^2\frac 1 {\e^2} (
      \theta^2+ 1)$$
      
      Since
      $$
      |\pi_N v|^2 = |v -\langle v , T\rangle T|^2= |v|^2 -(\langle
      v , T\rangle)^2$$
      then

  \begin{eqnarray*}
    |\pi_N\derpar v C_k|^2=
    \left| \frac 2 \e \theta \pi_N c_1 (\theta)  \right|^2
    =\frac 4 {\e^2} \theta^2 \Big( |c_1(\theta)|^2
      - (\langle T,c_1(\theta)\rangle)^2 \Big)
    =\\=
    \frac 4 {\e^2} \theta^2
    \left(1 - \langle\derpar\theta C,c_1(\theta)\rangle^2
      \frac 1{|\derpar\theta C|^2}
    \right) 
    =\\=
    \frac 4 {\e^2} \theta^2
    \left(1- \left\langle c_1(\theta)  2\frac v{\e},c_1(\theta)\right\rangle^2
      \frac 1{|\derpar\theta C|^2}   \right)
    =
    \frac 4 {\e^2} \theta^2
    \left(1 -  4\frac 1 {\e^2} v^2\frac 1{|\derpar\theta C|^2} \right)
    =\\=
    \frac 4 {\e^2} \theta^2\left(1 -  \frac 1{1+\theta^2} \right)
    =\frac 4 {\e^2} \theta^4 \frac 1 {1+\theta^2} 
  \end{eqnarray*}

  so that the energy for the part $v\in[0,1/2]$ of the homotopy is
  
  \begin{eqnarray*}
    E^N(C_k)= 4 k \int_0^{1/2} \int_0^\e
    \left| \pi_N \derpar v  C_k \right|^2
    \left|\derpar\theta C_k\right|~d\theta~dv  
    =\\ =
    4 k \int_0^{1/2} \int_0^\e
    \frac 4 {\e^2} \theta^4 \frac 1 {1+\theta^2}
    \left|\derpar\theta C_k\right|~d\theta~dv  
    =\\ =
    16 k \frac 1 {\e^2}    \int_0^{1/2}   \int_0^\e\theta^4
    \frac 1{1+\theta^2}
    \sqrt{4 v^2\frac 1 {\e^2} (   \theta^2+  1)}
    ~d\theta~dv   
    =\\=
    32 k \frac 1 {\e^3}    \int_0^{1/2}  \int_0^\e\theta^4
    \frac 1{1+\theta^2}
    v \sqrt{(   \theta^2+  1)}
    ~d\theta~dv
    =\\=
    32 k \frac 1 {\e^3} \frac 1 8   \int_0^\e\theta^4
    \frac 1{\sqrt{(   \theta^2+  1)}}   ~d\theta
    \le\\ \le 
    4 k \frac 1 {\e^3}  \frac{\e^5}5  = \frac 4 5 \pi^2  \frac 1 k
   \end{eqnarray*}
  
 \item 
   similarly in region $\theta\in[0,\e]$ $v\in[1/2,1]$ we have
   $$C_k(\theta,v) 
   = c_1(\theta) \left(1-2\frac{1-v}{\e} (\e -\theta)\right) 
   $$
   but we  implicitely  change variable $\theta\mapsto\theta-\e$,
   $v\mapsto v-1$
   to write
   $$C_k(\theta,v) 
   = c_1(\theta) \left(1-2\frac{v}{\e} \theta\right) 
   $$
   (this means that we will
   integrate on $\theta\in[-\e,0]$ $v\in[-1/2,0]$). Then

   $$\derpar v C_k = -c_1(\theta)2 \theta \frac 1 \e $$

   and

   $$\derpar \theta C_k 
   = \overdot c_1(\theta)\left(1-2\frac{v}{\e} \theta\right)
   - c_1(\theta)2\frac{v}{\e}
   $$
   
   and
   $$|\derpar \theta C_k|^2=\left(1-2\frac{v}{\e} \theta\right)^2
   +4\frac{v^2}{\e^2}
   =\frac{1}{\e^2} \Big((\e- 2v\theta)^2+4v^2\Big)$$

   \begin{eqnarray*}
     |\pi_N\derpar v C_k|^2=
     \left| \frac 2 \e \theta \pi_N c_1 (\theta)  \right|^2
     =\frac 4 {\e^2} \theta^2 \Big( |c_1(\theta)|^2
     - \langle T,c_1(\theta)\rangle^2 \Big)
     =\\=\frac 4 {\e^2} \theta^2 \left(1
       - \left\langle c_1(\theta)2\frac{v}{\e} 
         ,c_1(\theta)\right\rangle^2
       \frac 1 {|\derpar\theta C_k|^2}\right)
     =\\=\frac 4 {\e^2} \theta^2 \left(1
       - 4\frac{v^2}{\e^2}   
       \frac 1 {|\derpar\theta C_k|^2}\right)
     = \frac 4 {\e^2} \theta^2 \left(
       \frac {\left(\e- 2v\theta\right)^2}
     {(\e- 2v\theta)^2+4v^2} \right)
   \end{eqnarray*}

   Note that 
   \begin{equation}
     \label{eq: den pos}
     (\e- 2v\theta)^2+4v^2\ge \e^2(1+2v)^2+4v^2 \ge
     \frac{2\e^4}{(1+\e^2)^2}
   \end{equation}
   (the positive minimum is reached  at $v=-\e^2/(2+2\e^2), \theta=-\e$)

   Since 
   \begin{eqnarray*}
     |\pi_N\derpar v C_k|^2|\derpar \theta C_k|=
     \frac 4 {\e^2} \theta^2 \left(
       \frac {\left(\e- 2v\theta\right)^2}
       {(\e- 2v\theta)^2+4v^2} \right)
     \sqrt{\frac{1}{\e^2} \Big((\e- 2v\theta)^2+4v^2\Big)}
     =\\=
     \frac 4 {\e^3} \theta^2 \left(
       \frac {\left(\e- 2v\theta\right)^2}
       {\sqrt{(\e- 2v\theta)^2+4v^2}} \right)
     =\frac 4 \e \tau^2 \frac {\e^2(1- 2v\tau)^2}
     {\sqrt{ \e^2(1- 2v\tau)^2+4v^2}}
   \end{eqnarray*}
 
   where $\tau=\theta/\e$
   so the energy for the part $v\in[1/2,1]$ of the homotopy
   becomes

   \begin{eqnarray*}
     E^N(C_k)=2 k \int_{-1/2}^0\int_{-\e}^0
     |\pi_N\derpar v C_k|^2|\derpar \theta C_k|~d\theta~dv
     = \\ =
     2 k \int_{-1/2}^0\int_{-1}^0
     4\tau^2 
     \frac {\e^2(1- 2v\tau)^2}
     {\sqrt{ \e^2(1- 2v\tau)^2+4v^2}}
     ~d\tau~dv
     = \\ =8 \pi  
     \int_{-1/2}^0\int_{-1}^0
     \tau^2 
     \frac {(1- 2v\tau)^2}
     {\sqrt{ (1- 2v\tau)^2+4v^2/\e^2}}
     ~d\tau~dv
   \end{eqnarray*}

   since by \eqref{eq: den pos} the integrand is continuous and 
   positive, and it decreases pointwise when $\e\to 0$, then
   so $E^N(C_k)\to 0$ as $\e\to 0$ (by Beppo-Levi lemma,
   or Lebesgue theorem).

 \end{itemize}

\item 
as a second step, consider a generic smooth curve $c_1$; we could
approximate it with a piecewise smooth curve $c_1'$ , where
each piece in $C'$ is either contained in a sphere, or in a radious
exiting from $0$: by using the above homotopy on each spherical piece,
and translating and scaling each radial piece to $0$, we can
have a homotopy

\item then, given two generic smooth curves $c_0,c_1$, we can
  approximate them as above to obtain $c_0',c_1'$, and build
  a homotopy
  $$ c_0 \to c_0' \to 0 \to c_1' \to c_1$$
  this final homotopy can be built with small energy
\end{enumerate}
\end{Proposition}

\begin{Remark}
  \def\a{\alpha} \def\b{\beta} 
  More in general, if  $\a>\b>0$ , then the energy
  \[ \int_I |\pi_{N}\derpar v C |^\a |\dot C|^\b \]
  is ill-defined, as is shown in \ref{exa: E not lsc}.
\end{Remark}

\begin{Proposition}[$E$ is ill-posed]\label{inf E 0}
  \def\e{\epsilon}
  Consider the energy $E(C)$ associated to the metric 
  \eqref{eq: scalar quasi N}.
  Fix $c_0$ and   $c_1$ to be two regular curves.
  $\forall \e>0$, there is a homotopy $C\in \C$
  connecting $c_0$ to $c_1$ such that $E(C)<\e$.
  \begin{proof}
    Consider a homotopy defined as in the previous proposition; if we
    mollify it, 
    we can obtain a regular homotopy $C'$ such that
    $\|C'-C\|_{W^{1,3}}$ and $\|C'-C\|_{\infinity}$ are arbitrarily
    small; then we can reconnect $C'$ to $c_0$ to $c_1$ with a small
    cost, to create $C''$: by some direct
    computation, $E(C'')<2\e$.

    By using prop.~\ref{prop: C repar} on $C''$
    we obtain a $\widetilde C$ such that $E^N(\widetilde C)=E^N(C'')$, and
    since $\pi_{\widetilde N}\widetilde C=0$,
    $$E(\widetilde C)=E^N(\widetilde C)=E^N(C'')\le 2\e$$
  \end{proof}
\end{Proposition}

\section{Derivation of Flows}

In this section we show the details of the calculations of the
minimizing flows for both the $H^0$ and conformal energies.

\subsection{Some preliminary calculus}

First we develop in the following subsections some of the calculus
that we will need to work with the geometric parameters $s$
and $\v$ introduced in section \ref{sec: conformal}.

\subsubsection{Commutation of derivatives}

Note that the parameters $s$ and $\v$ do not form true coordinates
and therefore have a non-trivial commutator. The third parameter
$t$ will come into play later when we consider a time varying family of
homotopies $C(u,v,t)$ and take the resulting time derivative of either
the $H^0$ or the conformal energy along this family.

\begin{eqnarray}
  \pderiv{t}\,\pderiv{\v}
  &=& \nonumber
   \pderiv{t}
   \left(\pderiv{v}-\frac{C_u\cdot C_v}{C_u\cdot C_u}\,\pderiv{u}\right)
  \\ &=& \nonumber
   \pderiv{t}\,\pderiv{v}
  -\frac{C_u\cdot C_v}{C_u\cdot C_u}\,\pderiv{t}\,\pderiv{u}
  -\left(\frac{C_{ut}\cdot C_v + C_u\cdot C_{vt}}{C_u\cdot C_u}
         -2\frac{(C_u\cdot C_v)(C_{ut}\cdot C_u)}{(C_u\cdot C_u)^2}\right)
   \pderiv{u}
  \\ &=& \nonumber
   \left(\pderiv{v}-\frac{C_u\cdot C_v}{C_u\cdot C_u}\,\pderiv{u}\right)
   \pderiv{t}
  -\frac{
     C_u\cdot\left(C_{tv}-\frac{C_u\cdot C_v}{C_u\cdot C_u}C_{tu}\right)
    +C_{tu}\cdot\left(C_v-\frac{C_u\cdot C_v}{C_u\cdot C_u}C_u\right)
   }{C_u\cdot C_u}
   \pderiv{u}
  \\ &=& \label{eq:com_tv}
  \fbox{$\displaystyle
   \pderiv{\v}\,\pderiv{t}
  -\big(C_s\cdot C_{t\v}+C_{ts}\cdot C_{\v}\big)\pderiv{s}
  $}
  \\ \nonumber \\ \nonumber \\
  \pderiv{t}\,\pderiv{s}
  &=& \label{eq:com_ts}
   \pderiv{t}\left(\frac{1}{\|C_u\|}\,\pderiv{u}\right)
  =\frac{1}{\|C_u\|}\,\pderiv{t}\,\pderiv{u}
  -\frac{C_{ut}\cdot C_u}{\|C_u\|^3}\,\pderiv{u}
  =\fbox{$\displaystyle
   \pderiv{s}\,\pderiv{t}-C_{ts}\cdot C_s\pderiv{s}
  $}
  \\ \nonumber \\ \nonumber \\
  \pderiv{\v}\,\pderiv{s}
  &=& \nonumber
   \pderiv{v}\left(\frac{1}{\|C_u\|}\,\pderiv{u}\right)
  -C_v\cdot C_s \pderiv{s}\left(\frac{1}{\|C_u\|}\,\pderiv{u}\right)
  \\ &=& \nonumber
   \frac{1}{\|C_u\|}\,\pderiv{v}\,\pderiv{u}
  -\frac{C_{uv}\cdot C_u}{\|C_u\|^3}\,\pderiv{u}
  -C_v\cdot C_s \pderiv{s}\,\pderiv{s}
  \\ &=& \nonumber
   \pderiv{s}\,\pderiv{v}
  -C_{vs}\cdot C_s \pderiv{s}
  -C_v\cdot C_s \pderiv{s}\,\pderiv{s}
  \\ &=& \label{eq:com_vs}
   \pderiv{s}\left(\pderiv{v}-C_v\cdot C_s \pderiv{s}\right)
  +C_v\cdot C_{ss}\pderiv{s}
  =\fbox{$\displaystyle
   \pderiv{s}\,\pderiv{\v}+C_{\v}\cdot C_{ss}\pderiv{s}
  $}
\end{eqnarray}

\subsubsection{Some identities}

Here we write down some useful identities regarding various
derivatives of the homotopy with respect to the geometric
parameters $s$ and $\v$.

\begin{enumerate}
  \item $C_s\cdot C_s=1$
  \item $C_s\cdot C_{\v}=0$
  \item $C_{\v s}\cdot C_{\v}= C_{s\v}\cdot C_{\v}=-C_{\v\v}\cdot C_s$
  \item $C_{\v s}\cdot C_s=-C_{ss}\cdot C_{\v}$
  \item $C_{s\v}\cdot C_s=0$
  \item $C_{ss}\cdot C_s=0$
\end{enumerate}

\subsubsection{Commutation of derivatives with integrals}

Finally, we write down how to commute derivatives and integrals
when differentiating with respect to $t$ or $\v$.

\begin{eqnarray}
  \pderiv{t}\int_0^L f\,ds
  &=& \nonumber
   \pderiv{t}\int_0^1 f \|C_u\| \,du
  =\int_0^1 f_t \|C_u\|+ f(C_{ut}\cdot C_s)\, du
  =\int_0^L f_t + f(C_{ts}\cdot C_s)\, ds
  \\ &=& \label{eq:int_t}
  \fbox{$\displaystyle
   \int_0^L f_t - f_s(C_t\cdot C_s)- f(C_t\cdot C_{ss})\, ds
  $}
  \\ \nonumber \\ \nonumber \\
  \pderiv{\v}\int_0^L f\,ds
  &=& \nonumber
   \pderiv{v}\int_0^1 f \|C_u\| \,du
  =\int_0^1 f_v \|C_u\|+ f(C_{uv}\cdot C_s)\, du
  =\int_0^L f_v + f(C_{vs}\cdot C_s)\, ds
  \\ &=& \label{eq:int_v}
   \int_0^L f_v - f_s(C_v\cdot C_s)- f(C_v\cdot C_{ss})\, ds
  =\fbox{$\displaystyle
   \int_0^L f_{\v}- f(C_{\v}\cdot C_{ss})\, ds
  $}
\end{eqnarray}

\subsubsection{Intermediate Expressions}

The last step, before begining the flow calculation will be to
introduce a few ``intermediate'' expressions that will help keep the
the expressions in the upcoming derivations from becoming
too lengthy.

\begin{eqnarray}
  m        &=& C_{\v}\cdot C_{\v}
  \\
  m_t      &=& 2\,C_{\v t}\cdot C_{\v} = 2\,C_{t\v}\cdot C_{\v}
  \\
  m_s      &=& 2\,C_{\v s}\cdot C_{\v}= 2\,C_{s\v}\cdot C_{\v}
                                      =-2\,C_{\v\v}\cdot C_s
  \\
  m_{\v}   &=& 2\,C_{\v\v}\cdot C_{\v}
  \\
  m_{\v\v} &=& 2\,C_{\v\v\v}\cdot C_{\v}+2\,C_{\v\v}\cdot C_{\v\v}
\end{eqnarray}

\begin{eqnarray}
  M        &=& \int_0^L C_{\v}\cdot C_{\v} ds
  \\
  M_t      &=& \int_0^L 2\,C_{t\v}\cdot C_{\v}
                      + 2\,(C_{\v\v}\cdot C_s)(C_t\cdot C_s)
                      - m\,(C_t\cdot C_{ss})\,ds
  \\
  M_{\v}   &=& \int_0^L 2\,C_{\v\v}\cdot C_{\v} - m\,(C_{\v}\cdot C_{ss})\,ds
  \\
  M_{\v\v} &=& \int_0^L 2\,C_{\v\v\v}\cdot C_{\v}
                       +2\,C_{\v\v}\cdot C_{\v\v}
                       -m\,(C_{\v}\cdot C_{ss\v})
  \\ && \hspace{5mm} \nonumber
                       -m\,(C_{\v\v}\cdot C_{ss})
                       -4\,(C_{\v\v}\cdot C_{\v})(C_{\v}\cdot C_{ss})
                       +m\,(C_{\v}\cdot C_{ss})^2
               \,ds
  \\ \nonumber
           &=& \int_0^L 2\,C_{\v\v\v}\cdot C_{\v}
                       +2\,C_{\v\v}\cdot C_{\v\v}
                       +2\,(C_{s\v}\cdot C_{\v})^2
                       +m\,(C_{s\v}\cdot C_{s\v})
  \\ && \hspace{5mm} \nonumber
                       -m\,(C_{\v\v}\cdot C_{ss})
                       -4\,(C_{\v\v}\cdot C_{\v})(C_{\v}\cdot C_{ss})
               \,ds
\end{eqnarray}

\begin{eqnarray}
  L        &=& \int_0^L ds
  \\
  L_t      &=& \int_0^L -C_t\cdot C_{ss}\,ds
  \\
  L_{\v}   &=& \int_0^L -C_{\v}\cdot C_{ss}\,ds
  \\
  L_{\v\v} &=& \int_0^L -C_{\v\v}\cdot C_{ss} - C_{\v}\cdot C_{ss\v}
                       +(C_{\v}\cdot C_{ss})^2\,ds
  \\ \nonumber
           &=& \int_0^L C_{s\v}\cdot C_{s\v}\,ds - C_{\v\v}\cdot C_{ss}
\end{eqnarray}

\subsection{$H^0$ flow calculation}

We are now ready to begin the flow calculation. We'll start with the
case of the $H^0$ energy in this subsection and then proceed to
the conformal case in the followin subsection.

We begin by considering a time-varying family of homotopies
$C(u,v,t):[0,1]\times [0,1]\times (0,\infty)\to \real ^n$
and write the $H^0$ energy as
\begin{equation}
  E(t)=\int_0^1\int_0^L \big\|C_{\v}\big\|^2\,ds\,dv
      =\int_0^1 M dv
\end{equation}

Then the variation of $E$ is
\begin{eqnarray*}
  E'(t)
  &=& \int_0^1 M_t\,dv 
     =\int_0^1
        \int_0^L 2\,C_{t\v}\cdot C_{\v}
               + 2\,(C_{\v\v}\cdot C_s)(C_t\cdot C_s)
               - m\,(C_t\cdot C_{ss})\,ds
      \,dv \\
  &=& \int_0^1\int_0^L
        2\,\Big(C_{tv}-(C_v\cdot C_s)C_{ts}\Big)\cdot C_{\v}
      \,ds\,dv
     +\int_0^1\int_0^L
        C_t\cdot\Big(
          2\,(C_{\v\v}\cdot C_s)C_s-m\,C_{ss}
        \Big)
      \,ds\,dv \\
  &=& 2\int_0^1\int_0^1
        \,C_{tv}\cdot C_{\v}
      \|C_u\|\,du\,dv
    \\&&
     +2\int_0^1\int_0^L
        (C_v\cdot C_s)(C_t\cdot C_{\v s})
       +(C_{vs}\cdot C_s+C_v\cdot C_{ss})(C_t\cdot C_{\v})
      \,ds\,dv
    \\&&
     +\int_0^1\int_0^L
        C_t\cdot\Big(
          2\,(C_{\v\v}\cdot C_s)C_s-m\,C_{ss}
        \Big)
      \,ds\,dv \\
  &=& 2\int_0^1\int_0^1
       -\,(C_t\cdot C_{\v})(C_{uv}\cdot C_s)
       -\,C_t\cdot C_{\v v}\|C_u\|
      \,du\,dv
    \\&&
     +2\int_0^1\int_0^L
        (C_v\cdot C_s)(C_t\cdot C_{\v s})
       +(C_{vs}\cdot C_s+C_v\cdot C_{ss})(C_t\cdot C_{\v})
      \,ds\,dv
    \\&&
     +\int_0^1\int_0^L
        C_t\cdot\Big(
          2\,(C_{\v\v}\cdot C_s)C_s-m\,C_{ss}
        \Big)
      \,ds\,dv \\
  &=& 2\int_0^1\int_0^L
        C_t\cdot\Big(
         -C_{\v v}
         +(C_v\cdot C_s)C_{\v s}
         +(C_v\cdot C_{ss})C_{\v}
        \Big)
      \,ds\,dv
    \\&&
     +\int_0^1\int_0^L
        C_t\cdot\Big(
          2\,(C_{\v\v}\cdot C_s)C_s-m\,C_{ss}
        \Big)
      \,ds\,dv \\
  &=& -\int_0^1\int_0^L
        C_t\cdot\Big(
          2C_{\v\v}-2(C_{\v\v}\cdot C_s)C_s
         -2(C_{\v}\cdot C_{ss})C_{\v}
         +m\,C_{ss}
        \Big)
      \,ds\,dv
\end{eqnarray*}

In the planar case, $C_{\v}$ and $C_{ss}$ are linearly dependent
(as both are orthogonal to $C_s$) which means that
\begin{equation}
  (C_{\v}\cdot C_{ss})C_{\v}=(C_{\v}\cdot C_{\v})C_{ss}=mC_{ss}
\end{equation}
and therefore
\begin{equation}
E'(t)=
  -\int_0^1\int_0^L
    C_t\cdot\Big(
      2\big(C_{\v\v}-(C_{\v\v}\cdot C_s)C_s\big)-m\,C_{ss}
    \Big)
  \,ds\,dv
\end{equation}
by which we derive the minimization flow 
\[  C_t = 2\big(C_{\v\v}-(C_{\v\v}\cdot C_s)C_s\big) - m C_{ss}\]

\subsection{Conformal flow calculation}

We now define the $H^0_\phi$ energy as
\begin{equation}
  E_{\phi}(t)=\int_0^1 \phi(L) \int_0^L \big\|C_{\v}\big\|^2\,ds\,dv
      =\int_0^1 \phi M dv
\end{equation}
and compute its derivative as
\begin{eqnarray*}
  E'(t)
  &=& \int_0^1 \phi'L_t M+\phi M_t\,dv  \\
  &=& \int_0^1 \left(
        -\phi'M \int_0^L C_t\cdot C_{ss}\,ds
        +\phi
        \int_0^L 2\,C_{t\v}\cdot C_{\v}
               + 2\,(C_{\v\v}\cdot C_s)(C_t\cdot C_s)
               - m\,(C_t\cdot C_{ss})\,ds
      \right)\,dv \\
  &=& \int_0^1\int_0^L
        2\phi\Big(C_{tv}-(C_v\cdot C_s)C_{ts}\Big)\cdot C_{\v}
      \,ds\,dv
    \\&&
     +\int_0^1\int_0^L
        C_t\cdot\Big(
          2\phi(C_{\v\v}\cdot C_s)C_s-(\phi m+\phi'M)C_{ss}
        \Big)
      \,ds\,dv \\
  &=& 2\int_0^1\int_0^1
        \phi\,C_{tv}\cdot C_{\v}
      \|C_u\|\,du\,dv
    \\&&
     +2\phi\int_0^1\int_0^L
        (C_v\cdot C_s)(C_t\cdot C_{\v s})
       +(C_{vs}\cdot C_s+C_v\cdot C_{ss})(C_t\cdot C_{\v})
      \,ds\,dv
    \\&&
     +\int_0^1\int_0^L
        C_t\cdot\Big(
          2\phi(C_{\v\v}\cdot C_s)C_s-(\phi m+\phi'M)C_{ss}
        \Big)
      \,ds\,dv \\
  &=& 2\int_0^1\int_0^1
       -(\phi\,C_t\cdot C_{\v})(C_{uv}\cdot C_s)
       -C_t\cdot(\phi'L_v C_{\v}+\phi\,C_{\v v})\|C_u\|
      \,du\,dv
    \\&&
     +2\phi\int_0^1\int_0^L
        (C_v\cdot C_s)(C_t\cdot C_{\v s})
       +(C_{vs}\cdot C_s+C_v\cdot C_{ss})(C_t\cdot C_{\v})
      \,ds\,dv
    \\&&
     +\int_0^1\int_0^L
        C_t\cdot\Big(
          2\phi(C_{\v\v}\cdot C_s)C_s-(\phi m+\phi'M)C_{ss}
        \Big)
      \,ds\,dv \\
  &=& 2\int_0^1\int_0^L
        C_t\cdot\Big(
         -\phi'L_v C_{\v}-\phi\,C_{\v v}
         +\phi(C_v\cdot C_s)C_{\v s}
         +\phi(C_v\cdot C_{ss})C_{\v}
        \Big)
      \,ds\,dv
    \\&&
     +\int_0^1\int_0^L
        C_t\cdot\Big(
          2\phi(C_{\v\v}\cdot C_s)C_s-(\phi m+\phi'M)C_{ss}
        \Big)
      \,ds\,dv \\
  &=& -\int_0^1\int_0^L
        C_t\cdot\Big(
          2\phi'L_{\v} C_{\v}
         +2\phi C_{\v\v}-2\phi(C_{\v\v}\cdot C_s)C_s
         -2\phi(C_{\v}\cdot C_{ss})C_{\v}
         +(\phi m+\phi'M)C_{ss}
        \Big)
      \,ds\,dv
\end{eqnarray*}

In the planar case, $C_{\v}$ and $C_{ss}$ are linearly dependent
(as both are orthogonal to $C_s$) which means that
\begin{equation}
  (C_{\v}\cdot C_{ss})C_{\v}=(C_{\v}\cdot C_{\v})C_{ss}=mC_{ss}
\end{equation}
and therefore
\begin{equation}
E'(t)=
  -\int_0^1\int_0^L
    C_t\cdot\Big(
      2\phi\big(C_{\v\v}-(C_{\v\v}\cdot C_s)C_s\big)
     +2\phi'L_{\v} C_{\v}
     +(\phi'M-\phi m)C_{ss}
    \Big)
  \,ds\,dv
\end{equation}

which entails the flow
\[ C_t =  
      2\phi\big(C_{\v\v}-(C_{\v\v}\cdot C_s)C_s\big)
     +2\phi'L_{\v} C_{\v}
     +(\phi'M-\phi m)C_{ss}
 \]

\newpage

\tableofcontents 

\bigskip

\listoffigures

\bigskip

\providecommand{\bysame}{\leavevmode\hbox to3em{\hrulefill}\thinspace}
\providecommand{\MR}{\relax\ifhmode\unskip\space\fi MR }
\providecommand{\MRhref}[2]{%
  \href{http://www.ams.org/mathscinet-getitem?mr=#1}{#2}
}
\providecommand{\href}[2]{#2}

\end{document}